\input amstex
\documentstyle{amsppt}
\input epsf 
\input psfrag
\input amssym

\font\tencyr=wncyr10 
\font\sevencyr=wncyr7 
\font\fivecyr=wncyr5 
\newfam\cyrfam \textfont\cyrfam=\tencyr \scriptfont\cyrfam=\sevencyr
\scriptscriptfont\cyrfam=\fivecyr
\define\hexfam#1{\ifcase\number#1
  0\or 1\or 2\or 3\or 4\or 5\or 6\or 7 \or
  8\or 9\or A\or B\or C\or D\or E\or F\fi}
\mathchardef\Shch="0\hexfam\cyrfam 57

\define\RR{\Bbb R}

\define\HH{\Bbb H}
\define\SS{\Bbb S}

\define\ZZ{\Bbb Z}

\hyphenation{Liou-ville}

\NoBlackBoxes

\topmatter
\title  Parity of $\pmb{n}$-Frames with Application to non-Procrustean Orthogonalization
\endtitle
\rightheadtext {}
\author Jon A. Sjogren
\endauthor
\affil Towson University
\endaffil
\address Towson, Maryland
\endaddress
\date March 2020
\enddate
\dedicatory
QVALIS PATER TALIS FILIVS\\
William James Firey philosophiae doctor maxima cum laude MCMXXIII - MMIV\\
Brook MacDonald Firey philosophiae doctor magna cum laude MCMLI - MCMLXLIV
\enddedicatory
\endtopmatter
\document

\abovedisplayskip=13pt
\belowdisplayskip=13pt

Mainly we consider $n$-frames, that is, collections of a quantity `$n$' vectors in $\RR^n$.  The `non-compact space' consists of vectors $u_1,\ldots, u_n$ that are linearly independent but otherwise arbitrary. This space ${\Cal U}$ has the topology given as a subspace of $\left(\RR^n\right)^n$, hence possesses one or more path-components. The non-compact ${\Cal U}$ may be considered as homeomorphic to $Gl (n; \RR)$, matching each $u_j$ with the $j$-th column of the matrix representing $\vec{\mbox{\bf u}} = \left(u_1, \ldots, u_n\right)$. But as a rule we downplay group-structure in general (except for dimension 3, classical rotation groups) and matrix groups in particular.

A more careful exposition than ours would prove all the results by means of vector topology, without even requiring the concept of matrix multiplication with its algebraic properties. We do concede a ``Principle of Low Dimension''. This means that concepts that are avoided for higher dimensions, including determinants, homomorphisms of compact groups etc. are permitted in low dimensions. Thus we speak of the ``compact group'' of real orthogonal matrices $O(3)$,  and a possible ``homomorphism'' $\Phi: {\Cal H}_1 \to O^+(3)$ where ${\Cal H}_1$ represents the {\it compact} group of ``unit quaternions'', and $O^+(3)$ means ``principal path-component''. The mapping $\Phi$ should also be a ``topological covering'', probably the ``universal'' covering map.

So we forge ahead onto the ``compact space of frames'' consisting of vectors $v_1,\ldots, v_n$ where each $v_j$ has $\mbox{norm}  = 1$, and furthermore, the vectors are not only ``linearly independent'' as a set, but also ``pair-wise orthogonal'' so that the {\it dot product} $v_i \cdot v_j = \delta_{ij}$ (Kronecker delta).

Neglecting for the most part the fairly natural group structure, we denote nonetheless this compact space of frames by $O(n)$, otherwise known as the Stiefel manifold of orthonormal $n$-frames, $V_{n,n}$. We will eventually show (without using a general determinant or its continuity properties) that $O(n)$, $n \geq 1$, is not path-wise connected. In fact the topology of $O(n)$ and the non-compact ${\Cal U}(n)$, also written as $V_{n,n}^*$ are closely related.

Writing $O^*(n)$ as the Stiefel manifold $V_{n,n}^*$, see \cite{James}, permits consideration of other Stiefel manifolds such as $V_{n,k}^*$, the space of quantity $k$ independent vectors in Euclidean $n$-space $\RR^n$.

I. M. James exhibits a deformation retraction of non-compact $V_{n,n}^*$ to $V_{n,n} \simeq O(n)$, starting with the Gram-Schmidt formula which inputs some linearly independent collection $\mbox{\bf u} = \{u_1,\ldots, u_n\}$ and outputs an orthogonal frame $\mbox{\bf v} = \{v_1,\ldots, v_n\}$. A slight refinement of this method amounts to starting with a square, non-singular matrix $M \in Gl(n)$ and putting it into the form $M = Q R$, where $Q$ is an orthogonal matrix, and $R$ is non-singular and upper-triangular. A {\it path} is constructed from $R$ by deforming all of its (non-zero) diagonal elements, multiplying by a non-negative function $a_{ii}(t)$ of one  variable. We arrange that when $t=0$, $a_{ii}(t) = 1$ leaves the diagonal unaltered, but when $t$ reaches the value $1$, the modified $\tilde{r}_{ii}$ has been transformed either to $+1$ or to $-1$. Meanwhile, the upper-diagonal elements $r_{ij}, j > i$, have been sent continuously to $0$ by means of mulitplier-functions $a_{ij}(t)$. At the final step, $\overline{R}$ has $\{\pm 1\}$ values on the diagonal, so to create $\overline{R} = I$, we would start with $M = Q'R'$, where $Q'$ differs somewhat from $Q$ through a sequence of `reflections'. For details see a text of numerical linear algebra such as \cite{Golub \& van Loan}.

A more natural way to comprehend the deformation of $O^*(n)$ into (compact) $O(n)$ is to prove and in fact $O^*(n) = V_{n,n}^* \simeq V_{n,n} \times \RR^{n(n-1)/2}$, a homeomorphism. Then the $\RR^m$ factor, $m = \left(\matrix
   n+1 \\
   2
    \endmatrix\right)$ is a contractible (Euclidean) space and can be shrunk to $\{0\}$ (the origin), a singleton space. In this context we will not need to worry about straightening out a sequence of $+1$'s and $-1$'s.

To prove the Polar Decomposition in the real case, let $M$ be square of size $n \times n$ and non-singular, we want  $M = Q \exp T$, where $T$ is `self-adjoint' (real symmetric), and {\it both} $Q$ (orthogonal) and $T$ are uniquely determined from $M$. Furthermore, the entries of $Q$ and of $T$ depend continuously upon the values of the entries of $M$.

\head Spectral Theorem\endhead

To get started on Real Polar Decomposition, we examine important approaches to the finite-dimensional {\it spectral theorem}. This term is often used to mean that a {\it positive definite} real (Hermitean) operator can be diagonalized, with eigenvalues {\it positive}, by means of an {\it orthogonal} similarity.

Our main approach is to consider a self-adjoint (``real Hermitean'' or symmetric) square matrix and prove that it is orthogonally diagonalizable taking {\it real} values on the diagonal. In case the original matrix is {\it positive definite} (and hence {\it non-singular}) it will turn out that these real values are actually {\it positive}. This is the case needed when constructing the Polar decomposition for a real, square, invertible matrix.

It is well to establish that $A \in Gl(n)$ as a matrix is {\it self-adjoint} with respect to the Pythagorean inner product on $\RR^n$, exactly when it is symmetric.

Let $\mbox{\bf e}_{\text{1}} = (1, 0, \ldots, 0), \ldots, \mbox{\bf e}_{\text{n}} = (0, 0, \ldots, 1)$ be the standard (orthogonal) basis of $\RR^n$. Given the Pythagorean inner product, another orthonormal basis $\mbox{\bf f}_{\text{1}},\ldots, \mbox{\bf f}_{\text{n}}$ could be taken. The algebra of coefficients would look the same. Arbitrary vectors $\mbox{\bf x}$, $\mbox{\bf y}$ may be written
$$\align
\mbox{\bf x} &=  \sum_{i=1}^n \beta_i \mbox{\bf e}_{\pmb i} \\
\mbox{\bf y} &=  \sum_{i=1}^n \gamma_i \mbox{\bf e}_{\pmb i}
\endalign
$$
Then defining $\text{\mbox{\bf z}} = A\text{\mbox{\bf x}}$ we have
$$\mbox{\bf z}=  \sum_{i=1}^n \epsilon_i \mbox{\bf e}_{\pmb i}\qquad ,$$
then in fact $\epsilon_i = \sum a_{ik}\beta_k$ where $A = [a_{ij}]$. But for $A$ to be ``self-adjoint'' means (*) $(Ax, y) = (x, Ay)$. These expressions equal $\sum_{i,k=1}^{n,n} a_{ik} \beta_k \gamma_i$ and $\sum_{i,k=1}^{n,n} a_{ik}\beta_i\gamma_{k}$ respectively. In case always $a_{ik} = a_{ki}$, the quantities are equal.

Conversely if we fix $i, k$ and choose $\beta_j = \delta_{jk}$ with $\gamma_j = \delta_{ji}$ (Kronecker delta), the equation $a_{ik} = a_{ki}$ is obtained. Hence the ``symmetric'' (real Hermitean) and self-adjointness properties are {\it equivalent} for $A$. In the sequel a good linear algebra resource is \cite{Gel'fand}.

We next review two approaches to the Spectral Theorem that do not bring in the Fundamental Theorem of Algebra or its equivalents. Schur's decomposition could be  used (a matrix is congruent, and similar, to a nearly upper triangular matrix), but this result uses the general existence of a (complex) eigenvalue. As things stand, we construct any needed eigen-vectors and the like as we proceed. We prove that for $\|x\| =1$ ($x$ on the unit sphere $\SS^{n-1}$), the quadratic form $(Ax, x)$ assumes its minimum called $\lambda_1$ , at a vector called $w_{1}$, and $(\lambda_1,\ w_1)$ form an eigenvalue-eigenvector pair for $A$. For this a lemma is required.

\bigskip
\noindent
{\bf Lemma 1}\quad Given $B$ real self-adjoint, and also positive semi-definite, then if for a vector $e \in \RR^n$ there holds
$$(Be, e) = 0\qquad,$$
it follows that $Be = \vec{0}$, the zero vector.

\bigskip
\noindent
{\it Proof.} \quad See \cite{Gel'fand}, p. 127.\hfill $\blacksquare$

\bigskip
Next we find the pair $(\lambda_1, w_1)$. Since $\SS^{n-1}$ is compact, and $\theta (x):= (Ax, x)$ is continuous in $x$, $\theta (x)$ assumes a {\it minimum} value, say $\lambda_1$, at a unit vector $w_1$.

In that case $(Ax-\lambda_1 x, x) \geq 0$ for all $x \in \SS^{n-1}$. But the left-hand expression assumes the value  0 at $w_1$, hence by Lemma 1, also
$$(A-\lambda_1 I) w_1 = 0\qquad,$$
so indeed $w_1$ is valid as an eigen-vector for the minimal $\lambda_1$. Then notice that the subspace $V$ orthogonal to $w_1$ is also invariant under $A$. We may work with $A$ on $V \simeq \RR^{n-1}$ as before and obtain a minimal value $\lambda_2$ of $(Ax, x)$ for $x \in \SS(V)$, an $(n-1)$-sphere of  unit radius, together with its corresponding direction called $w_2$.

The properties of the original linear transformation now denoted by $B$, in the new basis $W = (w_1,\ldots, w_n)$ emerge without algebraic manipulation. Firstly, the eigenvalues are obviously {\it real} in the analytic formulation (find the maximum value attained by a quadratic form on a sphere). Secondly, the basis elements of $W$, scaled to unit norm, are pairwise orthogonal by construction.

``New coordinates'' expressed in the ``old coordinates'' (the {\it standard} coordinates) are given using matrix conjugation by $U = [w_1|w_n|\cdots|w_k]$, which yields an orthogonal matrix, so that
$$D = U^T A U$$
is just
$$D = \left[\matrix
\lambda_1 & \cdots & 0 \\
\vdots & \lambda_2 & \vdots \\
0 & \cdots & \lambda_n
\endmatrix\right]\qquad.$$
The only sophisticated aspect of the  above proof lies in the use of ``compactness'', the Bolzano-Weierstrass Theorem, applied to a high-dimensional space. In \cite{Geck}, the author shows how to find the minimal eigen-value (and corresponding eigen-vector) by using only the {\it completeness} of the real numbers, with an inequality
$$\|Av\| \leq n^{3/2} |A|_{\infty} \|v\|, \quad \mbox{for all $v$\qquad,}$$
which compares the Euclidean norm of vectors, and the {\it max norm} (of absolute entries $|a_{ij}|$).

We conclude this Section with a summary of the celebrated proof of the real spectral decomposition due to H. Wilf, enshrined in the archive ``Proofs from THE BOOK'' \cite{A-Z}. This proof does not mention eigen-values or eigen-vectors. Wilf does emphasize the ``matrix group''. Instead we will use the  space of $n$-frames as a collection of coordinate systems.

Recall the topology of the Stiefel manifold $V_{n,k} \subset \RR^n \times \RR^n \times \cdots \times \RR^n$ (with $k$ factors, $k \leq n$). The vectors of frame $\mbox{\bf f}$ in $V_{n,k}$ have unit norm, hence $V_{n,k}$ is bounded, also it is {\it closed} in the product topology, since the vectors are pair-wise orthogonal. A sequence of frames has a convergent sub-sequence, hence $V_{n,k}$ is {\it compact}, see \cite{James}.

Given a symmetric real matrix $A$, define $\Lambda(A)$ to equal the sum of the squares of its off-diagonal entries. Thus $\Lambda (A) \geq 0$. Now for each frame $\mbox{\bf f} \in V_{n,n}$, $A$ can be presented in
coordinates (entries) based on $\mbox{\bf f}$. Thus
$$A f_i = \sum_{i=1}^n a_{ij}^f f_j, \quad i = 1,\dotsc, n\qquad.$$
Here we may write explicitly
$$\Lambda (B_f) = \sum_{p, q\;p\ne q}^{n,n} a_{pq}^f \cdot a_{pq}^f\qquad,$$
where $B_f$ gives the realization of the linear transformation $A$ with respect to the given frame $\mbox{\bf f}$.

Fixing $A$, we may also write $H_A (\mbox{\bf f}) = \Lambda (B_{\text{\bf f}})$. This continuous function $H_A$ will attain a {\it real} minimum $\geq 0$ on the compact $V_{n,n} \simeq O(n)$. Suppose that the minimum is attained at the orthonormal frame $\mbox{\bf f}_{\text{0}}$. If $H_A (\mbox{\bf f}_{\text{0}}) = \rho = 0$, then the transformation $A$ in
$\mbox{\bf f}_{\text{0}}$-coordinates has the vectors constituting $\mbox{\bf f}_{\text{0}}$ as eigen-vectors. In other words, $B_{\text{\bf f}_{\text{0}}}$ is diagonal. The change-of-basis matrix $P_0$ has its columns the constituent vectors $\{v_0, v_1, \dotsc, v_n\}$ of $\mbox{\bf f}_{\text{0}}$, making $P_0$ into an orthogonal matriz and then
$$B_{\text{\bf f}_{\text{0}}} = P_0^T A P_0\qquad.$$

\noindent
A better notation might have been to write instead of $A$ in the above, $B_{\text{\bf e}}$,  where $\text{\bf e}$ is the standard coordinate-frame represented by the identity matrix $I_{n\times n}$.

In case $\rho = H_A (f_0) > 0$, the matrix $B_{f_0}$ would not be a diagonal matrix. But one may now find a new frame $g = f'$ such that $H_A(g) < \rho$, contradicting the minimality of $\rho$.

To avoid burdening the notation, take now $A$ (in place of $B_{\text{\bf f}_{\text{0}}}$) as a real symmetric matrix with $\Lambda (A) > 0$, so that $A$ is {\it not} a real diagonal matrix.

In particular suppose that $a_{rs} \ne 0$ where $r \ne s$. Construct an orthogonal matriz $U$, which looks like the identity $I$ except for the $2 \times 2$ sub-matrix at row and column indices $r$ and $s$. Hence, part of $U$ is of the form
$$ \left.\matrix
& r & s \\
r & \cos \theta & \sin \theta \\
s & -\sin \theta & \cos \theta
\endmatrix\right.\qquad .$$
Here $0 \leq \theta \leq \frac{\pi}{2}$ gives a parameter in radians.

Holding to our generic notation, we set $B = U^T A U$ and get a (two-sided) Jacobi ``rotation'' of $A$. Compare this with the Givens notation as developed further ahead.

A calculation now shows that for all $b_{ij}$, $i \ne j$, the values $b_{ij}^2 + b_{ji}^2 = a_{ij}^2 + a_{ji}^2$, except for $b_{rs}$ and $b_{sr}$. Thus only these two new entries contribute to a change in $\Lambda(A)$ toward its new value $\Lambda (B) = H_A (\text{\bf f}_{1})$, no matter what  the choice of the ``angle'' $\theta$. In fact, when $\theta$ is chosen $= 0$ radians, we obtain $b_{rs} = a_{rs}$, whereas for $\theta = \frac{\pi}{2}$, $b_{rs}$ will get the new value $-a_{rs} \ne 0$.

Therefore by the Intermediate Value Theorem applied to a continuous function $b_{rs}(\theta)$, the value ``zero'' must be
attained for some $\theta_0$ satistying $ 0 \leq \theta_0 \leq \frac{\pi}{2}$.

For this value it follows that $\Lambda \left(U_{\theta_0}^T A U_{\theta_0}\right) < \Lambda (A)$ is a strict inequality.
If we apply the above argument to the ``minimal'' (but non-diagonal) matrix $B_{f_0}$ as above, the given inequality immediately contradicts the minimality of $\rho = H_A(\text{\bf f}_{\text{0}})$.
Hence the {\it frame} $f_0$ where $\Lambda (B_{\text{\bf f}_{\text{0}}})$ attains its {\it minimum} ensures that  this minimum value $\rho$ must equal zero, and that $B_{\text{\bf f}_{\text{0}}} = P_0^T A P_0$ is a diagonal matrix.
This proof, consistent with our intentions in the remainder of the article, does not deal with the global group structure of $Gl(n; \RR)$ or $O(n)$, as done in \cite{A-Z}, but regards these objects as spaces of vector $n$-tuples, forming the columns of a matrix.

\head Polar Form and a Deformation Retract\endhead

We start with a real square matrix $A$ of order $n$. The case of interest is when $A$ is non-singular, hence we may consider $A \in Gl(n; \RR)$. We could write $A = RP$ where $R$ is orthogonal (it is square and its columns form an orthonormal set of vectors), with $P$  symmetric and positive semi-definite. The pair $(R, P)$ constitutes the {\it right} polar decomposition. When $A$ is invertible, so is $P$, giving $R$ as {\it positive definite}, with real eigen-values $\{\lambda_{\alpha}\}$, $\lambda_{\alpha} > 0$.

We collect some basic matrix facts.

\proclaim {Basic 1}
If $B$ is real symmetric with all {\it positive} eigen-values, $B$ is positive definite.
\endproclaim

\bigskip
\noindent
{\it Proof.} \quad
 By the Spectral Theorem from the previous Section, $B = U^T \Sigma U$ for a real {\it diagonal} $\Sigma$, and orthogonal $U$. We have

 \bigskip
 \hfill $X^T B X = X^T U^T \Sigma Ux = y^T \Sigma y = \sum_{i=1}^n \lambda_i y_i^2 > 0\qquad.\hfill \mbox{$\blacksquare$}$

\bigskip

\proclaim {Basic 2}
Next, if $B$ is real symmetric positive definite, then all of its eigenvalues are {\it positive}.
\endproclaim

\bigskip
\noindent
{\it Proof.} \quad If $(\lambda, X)$ is an eigen-value, eigen-vector pair for $B$,  we infer that wince $X$ is a non-zero vector, that $X^TB X = \lambda X^T X$. But $X^T B X > 0$, also $X^T X > 0$, hence $\lambda > 0$.\hfill \mbox{$\blacksquare$}

\bigskip
In the above assertions one may replace ``positive definite'' by ``positive semi-definite'' as long as one substitutes ``greater than or equal'' for ``greater than''.
The building  block of the Polar decomposition is the matrix $A^TA$. When $A$ is invertible, we have for any non-zero $v \in \RR^n$ that $Av$ is non-zero and $\langle v, A^T A v\rangle = \langle Av, Av\rangle > 0$.

Hence $A^T A$ is not only symmetric, but actually positive definite. It follows that $C = A^T A$ is non-singular. A different argument is that if $C$ were singular, then either $A$ would have a non-zero vector in its kernel, or else some non-zero $Au$ would be in the kernel of $A^T$, making the column rank of $A < n$. But row rank equals column rank, so $A$ was actually singular.

Conversely, if $A$ were non-invertible, evidently so is $C = A^T A$. Alternatively, simply take $B = A^{-1}\left(A^{-1}\right)^T$ and verify that $BC = I_{n\times n}$ using $\left(A^{-1}\right)^T A^T = \left(AA^{-1}\right)^T$, and that $CB = I_{n\times n}$ using $A^T \left(A^{-1}\right)^T = \left(A^{-1}A\right)^T$.

Now we may apply the Spectral Theorem with related discussion as above, to the $n\times n$ symmetric matrix $C = A^T A$. Since $C = V \Sigma V^T$, where $V$ is orthogonal and $\sigma$ is diagonal with real positive entries, we derive
$$\Sigma = \left[\matrix
\lambda_1 &  && 0 \\
 & \lambda_2 & &\vdots \\
 \vdots & &\ddots & \vdots \\
 &  &&\lambda_n
\endmatrix\right],\;\lambda_j > 0\qquad.$$
Then we may define a ``square root'' of $C$ by the formula
$$P:= \sqrt{C} = V\left[\matrix
\sqrt{\lambda_1} &  & 0 \\
 &   &\vdots \\
 \vdots & \ddots & \vdots \\
 &  &     \sqrt{\lambda_n}
\endmatrix\right] V^T\qquad.$$
Since the (lines generated by) the eigen-vectors of $P$ must have as eigen-values the respective square roots of the eigen-values of $C$. Therefore, $P$ is symmetric and positive definite, uniquely determined from $C$. Details of the construction can be found in Section 2.5 of \cite{Hall} and in \cite{Gavish}.

We put the previous work together in a

\bigskip
\noindent
{\bf Proposition 1 (Unique Polar Decomposition)}\quad
Any real invertible matrix $A$ may be written {\it uniquely} as
$$A = Re^X\qquad,$$
where $R \in O(n)$ (orthogonal), and $X \in \text{Symm}(n)$, symmetric but possible {\it singular}. The matrices $R, X$ (or their {\it entries}) depend continuously on this matrix $A$ (or its entries).

\bigskip
\noindent
{\it Proof.} \quad Simply define $R = AP^{-1} = A \cdot \left(A^T A\right)^{-1/2}$. We verify that $R$ is orthogonal, to wit
$$R^T R = \left(\sqrt{\left(A^T A\right)}^{-1}\right)^T A^T \cdot A \cdot \left(\sqrt{A^T A}\right)\qquad.$$
But we observed that both $A^TA$ and $\sqrt{A^TA}$ are self-adjoint. Hence the {\it inverse} $E$ of $\sqrt{A^T A}$ is also self-adjoint, we may now write
$$R^T R = \left(\sqrt{A^T A}^{-1}\right) \sqrt{A^T A} \cdot \sqrt{A^T A} \left(\sqrt{A^T A}\right)^{-1}\qquad,$$
which by the associativity of matrix multiplication must equal the identity $I_{n\times n}$. Hence we arrive at $A$ in the form of a product of orthogonal $R$ with positive definite $P$. To move from $P$ to symmetric $X$ we may use the local inverse of the matrix exponential $\exp(Y) = e^Y$:
$$\log Z = \sum_{m=1}^{\infty} \frac{(-1)^{m+1}}{m} \left(Z - I_{n\times n}\right)^m\qquad.$$
This series converges when $\|Z\| < \ln 2$ (real natural logarithm), where we indicate the Hilbert-Schmidt norm, `` the square root of the sum of the squares of the matrix entries of $Z$'', see \cite{Hall} Chpt. 2.

Furthermore, the logarithm of $P$ can be defined uniquely and continuously in $P$ by ``real analytic continuation''. The matrix $P$ can be {\it scaled} by a real factor to $P_0 = e^{-a}P$. The norm of $P_0$ should be chosen small enough for its logarithm to be uniquely defined and to vary continuously with its agrument. Now take
$$X = a I_{n\times n} + \log \left(e^{-a}P\right)\qquad.$$
Actually, $X$ ends up as {\it symmetric}, since given $Y$ symmetric, possibly singular, the matrix function $L = \exp (Y)$ is always defined, continuous, and produces a positive definite L.\hfill \mbox{$\blacksquare$}

\bigskip
From this Proposition follows immediately part of an important general ``diagonalization'' work-horse. We are referring to the Singular Value Decomposition in the ``real invertible'' case.

\bigskip
\noindent
{\bf Proposition 2}\quad
With notation as above, we obtain $A = W \Gamma V^T$, where $\Gamma$ is diagonal $\left(\Gamma = \sqrt{\Sigma}\right)$. The real matrices $W$ and $V$ are each {\it orthogonal}.

\bigskip
\noindent
{\it Proof.} \quad
We have by Proposition 1,
$$A = RP = RV \sqrt{\Sigma} V^T = W \sqrt{\Sigma} V^T\qquad,$$
where $W:= RV$.\hfill $\blacksquare$

\bigskip
We complete the topological decomposition of $Gl(n; \RR) \simeq V_{n,n}^*$.

\bigskip
\noindent
{\bf Theorem 1}\quad
The space of ortho-normal frames $V_{n,n}$ is a deformation rectraction of the (non-compact) space of linearly independent $n$-frames called $V_{n,n}^*$.

\bigskip
\noindent
{\it Proof.} \quad
In the previous Cartesian product, the factor $\RR^m$ is a vector space, hence contractible. The standard contraction of $\RR^m$ to $\{\vec{0}\}$ yields the desired deformation of $V_{n,n}^* \simeq V_{n,n} \times \RR^m$ to $V_{n,n} \simeq O(n)$.\hfill $\blacksquare$

\bigskip
We conclude that all homotopy properties of $V_{n,n}^*$ are the same as those of $V_{n,n}$. The two spaces $O^*(n)$ and $O(n)$ are ``homotopy equivalent''. Hence their respective sets of {\it path-components} are isomorphic as sets: they are of the same cardinality.

\head The Space of Orthonormal $n$-Frames\endhead

We observed in the previous Section that the (non-compact) topology of the linearly independent $n$-frames $V_{n,n}^*$ is essentially the same as that of the orthonormal $n$-frames  $V_{n,n}$ (according to homotopy equivalence). The latter space is conveniently represented, through a matrix or set of column vectors. The $n\times n$ real orthogonal matrices form a group $O(n)$, though for the most part we downplay the {\it product} and {\it inverse} operations on $O(n)$, (homeomorphic to $V_{n,n}$). Similarly on the linear algebra side, we tend to neglect the general theory of eigen-values and eigen-vectors, as well as ``determinants''. (Almost everything we deal with henceforth yields determinant $= \pm 1$.)

In particular we avoid the Fundamental Theorem of Algebra (or `Complex Axis Theorem', see \cite{Sjogren, Axis}), which is often used to concoct a needed real or complex eigen-vector.

We are interested in the path-connectedness of $O(n)$. For a space of its type: compact, Hausdorff and locally Euclidean of fixed dimension, namely $\left(\matrix
n \\
2
\endmatrix\right)$, the number of connected components in fact equals the number of path components, see \cite{Dugundji}.

We define $O^+(n)$ as the path-component in $O(n)$ of the $n\times n$ identity matrix, that is, the path-component containing the standard Euclidean frame {\bf e} with $e_i = (1, 0, \dotsc, 0), e_2, \dotsc, e_n$. At this point we wish to show that $V_{n,n} \simeq O(n)$ has at most {\it two} path-components. Later in the article, we employ the theory of covering spaces to confirm that $O(n)$ in fact possesses {\it exactly} two path-components. Hence  $O(n)$ contains $O^+(n)$ but is not the same space as $O^+(n)$.

Given a frame $\text{\bf v} \in V_{n,n}$ represented through {\it columns} as an orthogonal matrix $M$, we connect {\bf v} with a {\it path} either to the ``standard frame'' $\text{\bf e} = \{e_1, e_2, \dotsc, e_n\} \sim I$ or to the ``toggled frame'' $\text{\bf e}^- = \{e_1, e_2, \dotsc, -e_n\}\sim I^{-}$. The two (related) methods to construct such a path are
\roster
\item"(i)"
a geometric approach from ex. 13, p. 28 of \cite{Hall}, and
\item"(ii)" the precise ``numerical'' approach of Givens rotations.
\endroster

As in Chpt. 1 of \cite{Hall}, when $n \geq 2$ there exists a continuous concatenation of infinitesimal ``rotations'' of $O(n)$ that lead from the first column $v_1$ of {\bf v} to the standard vector $e_1$. Suppose firstly that $v_1$ happens to equal one of the standard $n$-frame vectors $e_j$. Then re-order the vectors of frame {\bf v} so that $v_1$ comes first as a column of the modified matrix $\hat{M}$. In case $v_1 = e_1$, no rotation is necessary to bring $v_1$ to $e_1$.

If $j \ne 1$, a quarter-turn of $90^{\circ}$ in the plane of $P_j = (e_1, e_j)$ is sufficient to take $v_1 = e_j$ to $e_1$. This ``path'' can be expressed through pre-multiplication by a continuum of Givens matrices
$$U_{\theta} = \left.\matrix
& 1 & j \\
1 & \cos \theta & \sin \theta \\
j & -\sin\theta & \cos \theta
\endmatrix\right. \qquad,$$
where $\theta$ ranges from $0$  to $\frac{\pi}{2}$ radians. Hence $U_{\theta}{M}$ now has $e_1$ as its first column.

More generally, if $v_1$ is not equal to any $e_j$, consider the plane $Q_1 = (v_1, e_1)$ in $W = \RR^n$. Starting with the (non-orthonormal) basis $\{e_1, v_1, e_3, e_4, \dotsc, e_n\}$ we use a modified Gram-Schmidt procedure to construct a new $n$-frame {\bf f} with $f_1 = e_1, f_2, \dotsc, f_n$ where the plane $(f_1, f_2)$ is the same as $(e_1, v_1)$. See \cite{Golub \& van Loan}. Without loss of generality, we may take $v_1$ in the first quadrant of $Q_1$. Now there is a Givens path with non-diagonal entries only from
$$\align
\text{row} &= \text{1  or  2} \\
\text{column} &= \text{1 or  2}\qquad, \\
\endalign
  $$
which is the {\it identity} at $\theta = 0$ and moves $f_2$ to $f_1 = e_1$ at $\theta = \frac{\pi}{2}$. Hence by continuity there exists $0 \leq \theta_0\leq \frac{\pi}{2}$ where the rotation $U(\theta_0)$ takes $v_1$ to $e_1$.

{\it Note:}\quad With terminology from topology, we may consider the planar rotation $U(\theta)$ as extending to a mapping
$$S^{n-2} U: \RR^n \to \RR^n$$
giving the ``$(n-2)$-fold  suspension'' of $U(\theta)$.

In any case, for the new coordinates, $\hat{M} = U(\theta_0) M$ has entries zero in the first row and column, except for $\hat{M}_{11}=1$.
The residual $M' \in O(n-1)$ acts on $\RR^{n-1} = \text{span}\{f_2,\dotsc, f_m\}$. The argument above can be repeated upon $M'$, and it breaks down only for $W_1 = \RR^1$, where the vector $-e_1$ cannot be rotated into $e_1$. But the
original matrix $M$ can be rotated by
$$U_{n-1}\left(\theta^{n-1}_0\right) U_{n-2} \left(\theta^{n-2}_0\right) \dotsc U_1\left(\theta_0^1\right)$$
either into $I_{n\times n}$ or into
$$I'_{n\times n} = \left.\matrix
1 & 0 & \cdots & 0 \\
0 & 1 & \\
\vdots & &\ddots & \\
0 & &&-1 &
\endmatrix\right.\qquad,$$

\noindent
representing the frame $e_{1},\dotsc, -e_n$. These two frames give both path components of $O(n)$, but we have not quite proved that there exists between these representatives {\it no path}. Hence for the moment we must admit the possibility that $O^+(n) = O(n)$. The first is often called ``special orthogonal group''.

\head Numerical Algebra Viewpoint\endhead

Standard texts such as \cite{Y-G} emphasize ``Givens rotations'' through their entries (generally four non-zero
entries), rather than by the {\it angle} of the rotation.

For our purposes a Givens rotations {\it really is} a continuous orthogonal transformation, as the rigid motion of a ball with fixed center. A ``Givens loop'' will be the critical algebraic item with which to understand the equivalence of $n$-frames.

The numerical process of zeroing entries of a matrix by means of rotations often gives way to the ``Householder transformation'' which is essentially a reflection about a hyperplane in $\RR^n$. But we saw that rotating $v_1$ to $e_1$ is just as efficient, at least within the class of orthogonal matrices ($n$-frames).


Writing the Givens transformation through angle $\theta$ now as $G(\theta)$, we have
$H(\theta) = G_{n-1}(\theta) \circ \cdots\circ G_1(\theta)$, which should be considered
as the concatenation of paths in $O(n)$, a {\it product} in the sense of ``Brandt Groupoid'' (the endpoint of $G_2 (\theta^2)$ should be the initial point of $G_3(\theta^3)$  and so forth). The total path $\gamma (\tau)$ might not be {\it smooth} in the variable $\tau$, but it is of course {\it continuous}.

The construction of the indices and delimiting angles of the segments $G_1(\theta), G_2(\theta),\dotsc$ is a matter of standard computational matrix theory. See ``Lectures on Linear Algebra'', \cite{Baker}, Lecture 8. ``Given'' two non-zero coordinates of $v_1$, at ``place'' $i < j$, we write
$$\left(\matrix
c & s \\
-s & c
\endmatrix\right) \left.\matrix
v_{1i}\\
v_{1j}
\endmatrix\right. =  \left(\matrix
\rho \\
0
\endmatrix\right) \qquad.$$
The preferred value of $\rho$   is a subject of numerical nit-picking, but it can be chosen as
$$\rho = \text{sign}\left(\max \left(v_{1i}, v_{1j}\right)\right) \cdot \left\|\matrix
v_{1i}\\
v_{1j}
\endmatrix\right\|\qquad.$$
When choosing the actual rotation $G\left(\theta_0'\right)$, it is important to carry out arithmetic algorithms that avoid digital underflow and overflow as may result from a ``na\"{\i}ve application''
$$\rho = \left(v_{1i}^2 + v_{1j}^2\right)^{\frac{1}{2}},\quad c = v_{1i}/\rho, \quad s = -v_{1j}/\rho$$
Of importance as we consider continuous Givens paths, the continuity of the procedure was emphasized by E. Anderson in ``A Working Note on LAPACK Revision'' (University of Tennessee, 2000). The improved Givens programs have numerically demonstrated their robustness.

The study of topology on a (``good'') space is certainly more convenient when the space is {\it path-connected}. At this point we are not yet sure whether $V_{n,n}$ is path-connected, so we work with its (path-) connected component $O^+(n)$. Next we give a constructive characterization of this {\it compact} space.

The numerical analysts have constructed Givens rotations and, unwittingly, Givens paths and Givens loops without a need for modification of the coordinate system. ``Given'' the orthogonal matrix $M$ with initial column $v_1$, we find a finite sequence of Givens paths so that $\hat{M} = G_{n-1}\left(\theta^{n-1}_0\right)\cdots G_1\left(\theta_0^1\right) \cdot M$, the new first column is all zeros, save for a single `$1$' in the $j$-th place. Our final canonical path
$H\left(\frac{\pi}{2}\right)$ effects a rotation taking $e_n$ to $e_1$. Then the resulting matrix $\check{M} = H(\theta) {M}$ has in its first column $e_1 = (1, 0, 0, \dotsc, 0)^T$.

\head Explicit Rendering of $O^+(n)$\endhead

A notation similar to what we use is found in \cite{Whitehead, Elements}. As also pointed out in \cite{James}, the path-component of $\text{\bf e} = \{e_1, e_2,\dotsc, e_n\}$ in $O(n)$ is homeomorphic to $V_{n, n-1}$.

The latter space does not afford a multiplicative (group) structure equally\break obvious as the one on $V_{n,n}$. We secretly admit that $O^+(n)$ is the same as ``special orthogonal group'', the group-structure somehow does reside there.

Given $M = \left(v_1^T, v_2^T,\dotsc, v_n^T\right)$, let $\Phi(M)$ be the $n\times (n-1)$ matrix
$$N = \left(v_1^T, v_2^T, \dotsc, v_{n-1}^T\right)\qquad.$$
We have simply dropped the final column vector from $M$. Notice that $\Phi:  V_{n,n} \to V_{n,n-1}$ must be some kind of fibration. Actually, $\Phi$ is a covering map: the inverse image $\Phi^{-1}(p)$ of a point consists of finitely many points, that is, finitely many $n$-frames. The covering is a {\it double cover} since there are exactly two unit vector solutions $v_n$ satisfying $\Phi(\pmb{q}, v_n) = \pmb{q}$, given $\pmb{q} \in V_{n, n-1}$.

One way to compute $v_n$ is by means of the ``general Massey cross-product'' using determinants in a manner similar to the Laplace expansion \cite{Massey, Cross}. Whatever answer $v_n$ that is obtained, its negative the vector $-v_n$ will serve equally well.

Without loss of generality one may take the vectors $v_1,\dotsc, v_{n-1}$ as spanning some linear hyperplane $\RR^{n-1} \subset \RR^n$.
Then any completion to $\left(v_1,\dotsc, U\right) \in V_{n,n}$ must have $U$ either in the upper- or lower-hemisphere
$$H_n^+ = \left\{u \in \RR^n\Big| u \cdot u = 1,\quad u_n > 0\right\}$$
or
$$H_n^- = \left\{u \in \RR^n\Big| u \cdot u = 1,\quad u_n < 0\right\}\qquad.$$

\newpage
\bigskip
\centerline{
\psfrag{an}{$\alpha \Cal N$}
\psfrag{u}{$\Cal U$}
\psfrag{p}{$P$}
\psfrag{n}{$\Cal N$}
\psfrag{b}{$b$}
\psfrag{O}{$O$}
\psfrag{R}{$\RR^{n-1}$}
\psfrag{e1}{$\alpha \Cal N + b = \Cal U$}
\epsfbox{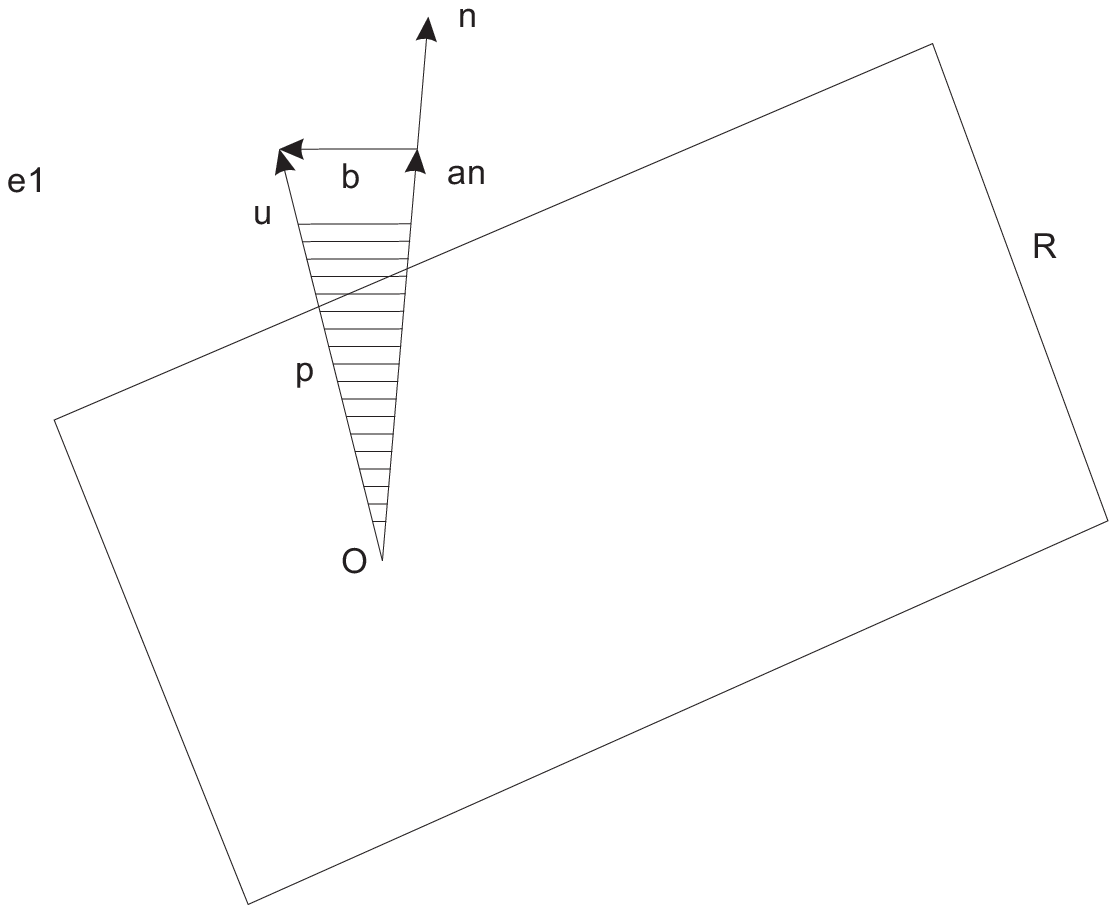}
}

\medskip
\centerline{Figure 1}

\bigskip

The non-zero dot product value contradicts the assumption that {\bf U} is perpendicular to all of $\Bbb R^{n-1}$.

Supposing that $\Cal U \in H_n^+$ but $\Cal U \ne \Cal N$ where $\Cal N = (0,0,\dotsc,1)$, then we use a Euclidean geometry argument in the plane $P \subset \RR^n$ that contains vectors $\Cal U$ and $\Cal N$. The angle $\epsilon = \angle (\Cal U,\Cal N)$ cannot equal $0$, so if $\vec{b}$ is a vector in the plane $P$ which is orthogonal to $\Cal N$, $\Cal U \cdot b = |\Cal U| \cdot |b| \cos \delta$, where $\delta$ is the {\it supplement} to $\epsilon$. Thus $\Cal U \cdot b \ne 0$ which contradicts the construction of $\Cal U$.

Thus we have

\bigskip
\noindent
{\bf Proposition A}\quad
The projection $\Phi: V_{n,n} \to V_{n,n-1}$ is a two-to-one covering map. In particular, $\Phi$ is locally a homeomorphism.

\bigskip
\noindent{\it proof.}\quad The latter property can be expressed in several ways, that given $x \in V_{n,n-1}$ and $\{y_1, y_2\} = \Phi^{-1}(x)$, with $y_1 \ne y_2$, any metrically small open neighbourhood  $N$ of $x$ yields $\Phi^{-1}(N) = R_1 \cup R_2 \subset V_{n,n}$.

Here $y_1 \in R_1, y_2 \in R_2$, $R_1$ and $R_2$ are disjoint, both $\Phi: R_1 \to N$ and $\Phi: R_2 \to N$ are homeomorphisms     (``Stack of Hot-Cakes'', see \cite{G-P}). See Figure 2.\hfill $\blacksquare$

\newpage
\centerline{
\psfrag{1}{Case 1}
\psfrag{2}{Case 2}
\psfrag{3}{$\Phi$ $\matrix
\text{double} \\ \text{covering}
\endmatrix$}
\psfrag{4}{$V_{n,n-1}$}
\psfrag{5}{$\Phi$}
\epsfbox{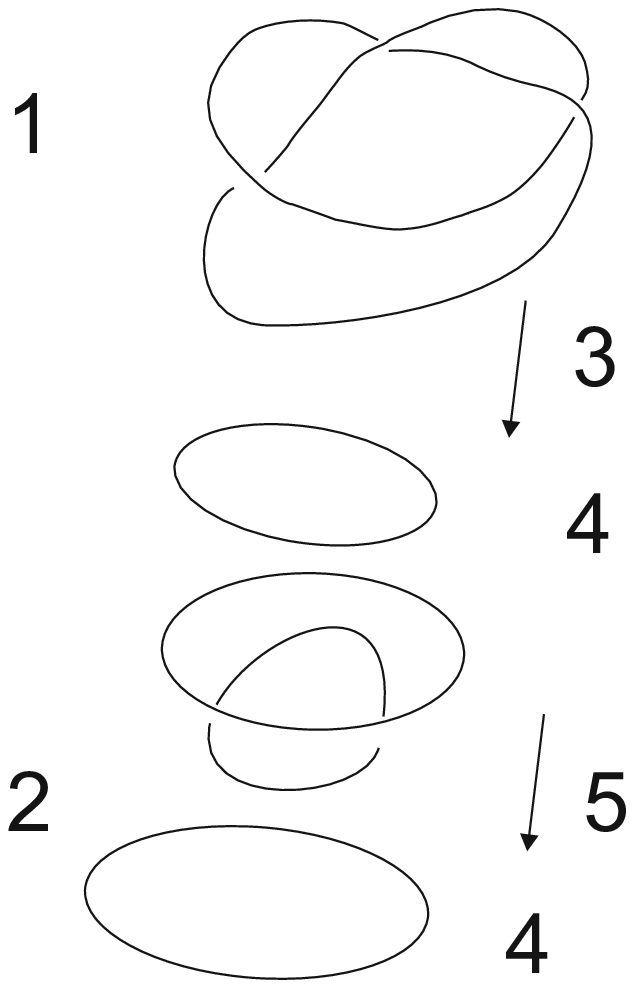}
}

\medskip
\centerline{Figure 2}

\bigskip
\noindent
{\bf Proposition B}\quad
The $(n,n-1)$ Stiefel manifold $V_{n,n-1}$ is path-connected for $n \geq 3$.

\medskip
\noindent
{\it Proof.} \quad
To show that $O(n)$ has at most 2 path-components, we used a Givens argument dealing with rotations in a plane generated by two special vectors of $\RR^n$. To treat $V_{n,n-1}$ we need extra room, so rotations in 3-dimensional subspaces come up. By previous arguments, it suffices to show that the two-frame $(e_1, e_2)$ in 3-space is rotatable to $(e_2, e_1)$. Indeed, one may obtain the needed `rotation' of $\RR^n$ by extending the 3-D rotation to $\RR^n$ by means of an $(n-3)$-fold suspension transformation. But $(e_1, e_2, e_3)$ is rotatable to $(e_1, -e_2, -e_3)$, which has a path to $(e_2, e_1, -e_3)$. But in $V_{3,2}$ this is the {\it same point} as $(e_2, e_1, e_3)$.

Hence the separate 2-frames have a path running between them in $V_{3,2}$.\hfill $\blacksquare$

\bigskip
We may conclude also that for $n \geq 3$, if $O(n)$ is {\it not} path-connected $(O(n)$ has two distinct path-components), then $V_{n,n-1} = O^+(n)$, that is, the path-component containing the canonical $n$-frame $\text{\bf e} = (e_1, e_2, \dotsc, e_n)$.

\bigskip
\noindent
{\bf Remark}\quad
``How we always knew that $\pi_0 (O) = {\Shch} (\ZZ_2)$, the set of two elements.''$^{*}$\vfootnote{$^{*}$}{Serge\u{\i} Olyegavich Shchakr'yennov ``On abstract Verdier complexes''. Plenary Session, Mongolian Academy of Sciences, Urga 1978.}\pagebreak

\head The Three Coverings\endhead

Assuming that we have at hand a sufficient understanding of $O(2)$, we may consider $O(3)$ to be the foundational case, working up from there to study $O(n)$, $n > 3$. We do not make significant use of the determinant as a {\it function}, the Theorem on Complex Polynomials (FTA), Lie groups and their quotients, or the {\it local group} (Lie algebra). We stay away from {\it simplicial, cubical} or {\it singular} homology theories, that afford a built-in concept of orientation.

What we know about $O(2)$, also known as the compact Stiefel space $V_{2,2}$, is that it has two path components, with the principal path-component (containing the canonical 2-frame) homeomorphic to a {\it circle}, so that $\pi_1(O^+(2)) = \pi_1 (V_{2,1}) \simeq \ZZ$, the infinite cycle group. We know that the universal covering space of $V_{2,1}$ is then the real line $\RR^1$, the case where this space is {\it not} compact. In fact, the only other universal cover that we examine is $\Phi: \SS^3 \to O^+(3)$ which enables us to compute $\pi_1 (O^+(3)) = \pi_1 (V_{3,2})$, the ``connectivity of loops on the space'', up to homotopy. Calculation of this fundamental-group will also enable us to determine $\pi_1(O^+(n)) = \pi_1 (V_{n, n-1})$ using the exact homotopy sequence of a particular fiber bundle (a small part, giving a {\it short}  exact sequence). The result of all these $\pi_1$ reckonings is the same: the ``toggle group''
of two elements, $\ZZ_2$.

The universal cover of $V_{3,2}$ will be examined in the next few Sections. First we will look at the pseudo-covering space

$$\left.\matrix
& V_{n,n} \\
& \rho_n \downarrow \\
& V_{n,n-1}, \\
\endmatrix\right.\qquad$$
so-called as the covering space $V_{n,n}$ is not (necessarily) path-connected.
But this is a double covering. A specific generator $\sigma$ of $\pi_1 (V_{n, n-1})$, $n>3$, can be constructed geometrically; in fact it can be represented by a Givens loop. Using suspension, the action of $\sigma$ on the space of $n$-frames is ``the identity''. Hence by a basic result in
(pseudo-) covering space theory, the covering $V_{n,n}$ {\it cannot} be path-connected.

Now we understand why the 4-frame $(e_1, e_2, e_3, e_4)$ in $V_{4,4}$ cannot be moved continuously, within  the orthonormal frames, to the ``flipped'' 4-frame $(e_2, e_1, e_3, e_4)$. The double covering $\rho_n = V_{n,n} \to V_{n, n-1}$ was mentioned earlier and is defined by erasing the final $n$-vector of the given frame. Thus $\rho_n$ applied to
$$\left.\matrix
| & | & | & | \\
v_1 & v_2 & v_3 & v_4 \\
| & | & | & | \\
\endmatrix\right.$$
in $V_{4,4}$ yields
$$\left.\matrix
| & | & |  \\
v_1 & v_2 & v_3  \\
| & | & |  \\
\endmatrix\right.$$
considered as an element of $V_{4,3}$. We saw that the spaces $V_{n+1,n}$ $n \geq 2$, are path-connected, hence suitable for a r\^ole as fiber, or one of {\it total space}, for a fiber bundle.

The third ``covering'' (in a loose sense) that we require is this fibration or ``fiber bundle'' structure. Explicitly, define $p^n : V_{n,n-1} \to \SS^{n-1}$ by taking the {\it initial} vector of the $(n-1)$-frame $\mbox{\bf v} = \left(v_{1},\dotsc, v_{n-1}\right)$, so that $p^n{\mbox{\bf v}} = v_1$, which is a unit vector of $\RR^n$, hence lies in $\SS^{n-1}$.

We follow Sect. 4.2 of \cite{Hatcher} concerning the exact homotopy sequence.  Another treatment that emphasizes the ``fibration'' and its consequences for the topology of Stiefel manifolds, can be found in IV.9-11 of \cite{Whitehead, Elements}.

Our mapping $p^n : V_{n,n-1} \to \SS^{n-1}$ gives rise to a fiber bundle. The inverse image of a metrically small neighborhood $Y \subset \SS^{n-1}$ containing $s \in \SS^{n-1}$ consists of $(s', u_2, \dotsc, u_{n-1})$ where $s' \in Y$ and $u_2, \dotsc, u_{n-1}$  together form an $(n-2)$-frame lying in the {\it equator} $E_{s'}$ of $\SS^{n-1}$ of unit vectors orthogonal to $s'$. For a given $s'$, this space $\{\mbox{\bf u'} = (u_2, \dotsc, u_{n-1})\}$ is linearly homeomorphic to $V_{n-1, n-2}$, so $\left(p^n\right)^{-1} (Y) \simeq Y \times V_{n-1,n-2}$.

We have at hand $V_{n,n-1}$ as a fiber bundle projecting to $\SS^{n-1}$, with $V_{n-1,n-2} \simeq O^+(n-1)$ as the fiber.

\bigskip
\noindent
{\bf Definition}\quad Given a general fiber bundle
$$(F, x_0) \hookrightarrow (X, x_0) \buildrel\hbox{$p$}\over{\longrightarrow} (B, b_0)$$
of {\it nice} spaces where base points  are indicated, there exists a {\it long exact} sequence of (absolute) homotopy groups (for which we display the lower terms)
$$\align
\cdots & \longrightarrow \pi_2 (B) \longrightarrow \pi_1 (F) \longrightarrow \pi_1 (X) \\
& \longrightarrow \pi_1 (B) \longrightarrow \pi_0 (F) \longrightarrow \pi_0 (X) \longrightarrow 0\qquad.\\
\endalign$$
To save space, the base points have been suppressed.

We know all about $\pi_1(X)$ when $X = O^+(3) = V_{3,2}$, so the  interesting cases are where $n \geq 3$ and $F = V_{n, n-1}$, $X = V_{n+1,n}$ and $B = \SS^n$.
All of these spaces are path-connected $B = \SS^n$ is simply-connected and most likely ``aspherical'' with respect to mappings $\tau: \SS^2 \to B$ of the {\it classical} sphere $\SS^2$. In case we are not sure about the evaluation of $\pi_2(B)$ we write for this ``weak case'' of the sequence
$$\align
  \longrightarrow &\pi_2 (\SS^n)  \buildrel\hbox{$\partial$}\over{\longrightarrow}\pi_1 (V_{n,n-1})  \buildrel\hbox{$\iota_{*}$}\over{\longrightarrow} \pi_1 (V_{n+1,n}) \\
 & \buildrel\hbox{$p^n_{*}$}\over{\longrightarrow} \pi_1 (\SS^n) = 0\qquad.
\endalign$$

\bigskip
\noindent
{\bf Proposition C}\quad
This ``weak assumption'' is sufficient to show that the orthogonal group $O(n+1)$ has exactly two path-components!
Thus path-distinct right-hand and left-hand $(n+1)$-frames do exist, say $(e_1,\dotsc,$
$ e_{n+1})$ and $(-e_1,\dotsc, e_{n+1})$.

\bigskip
\noindent
{\it Proof.} \quad
We will compute $\pi_1 (V_{3,2}) \simeq \ZZ_2$ in the following Section. After that, assume by an induction hypothesis that $\pi_1 (V_{n, n-1})$ is either $\ZZ_2$ or the $0$ abelian group. Then whatever $\pi_2 (\SS^n)$ may be (but {\it abelian}), except for the $0$ group, $\pi_1 (V_{n+1, n})$ turns out to be $0$. Otherwise, when $\pi_2(\SS^n) = 0$ as expected, we have the isomorphism
$$\pi_1 \left(O^+(n)\right) \buildrel\hbox{$\iota_{*}$}\over{\longrightarrow}\pi_1 \left(O^+(n+1)\right)$$
as an isomorphism of $\ZZ_2$-groups, induced by the inclusion mapping
$$\left.\matrix
v_1 & v_2 & \cdots & v_{n-1} &\mapsto &  v_1 & v_2 & \cdots & v_{n-1} & 0\\
| & | &  & | & &  | & | &  & | & | \\
  &   &  &   & &  0 & 0 &  & 0 & 1 & .
\endmatrix\right.$$
Later, we construct a compatible set of generators for $\left\{\pi_1\left(O^+ (n)\right)\right\}$. In fact as a path $\gamma (\theta)$ in $O^+(n)$, the {\it induced path} $\gamma^*(\theta)$ in $O^+(n+1)$ is just the same path $\gamma$, keeping in mind the above inclusion $\iota : O(n) \hookrightarrow O(n+1).$

In case $\pi_1 \left(O^+(n+1)\right) = 0$, then $O(n+1)$ is {\it not} connected, since in the  double covering $\rho = O (n+1) \to O^+(n+1)$, there is no action of the base loops upon the discrete set of lifted ``base points'' $\kappa = \left\{\rho^{-1}(x_0)\right\}$: the loops act as the identity. If $O(n+1)$ {\it were} path-connected then each $y \in \kappa$ could be transported to its double $y' \in \kappa$ by means of a base loop action. See \cite{Munkres}, (Topology) Theorem 54.4.\hfill $\blacksquare$

\bigskip
In the remaining case where all $\pi_1(V_{k , k-1})$  equal $\ZZ_2$, $k = 3,\dotsc, n+1$, we use a similar argument in terms of a canonical generator of $\ZZ_2$. The argument will  be recapitulated in the next Section, along with the determination of $\pi_1 \left(O^+ (3)\right)$, the homotopy classes of loops on the traditional ``rotation group''. An explicit geometrical homeomorphism is often given between this space and the ``projective 3-space'' $\RR P^3$, but we avoid a construction which has fairly been beaten to death.  Instead we use computational algebra, relying for the Rectitude of our Actions on the Principle of Low Dimension, which lets us ``get away'' with cross-product, groups and homomorphisms, concepts that we prefer not to use when the ``manifold'' dimension number is unbounded.

Now there is an advantage toward settling on $\pi_2 (\SS^2) = 0$, ``the zero group'', for $n \geq 0$ together with the isomorphism
$$0 \to \pi_1(V_{n,n-1}) \buildrel\hbox{$\iota_{*}$}\over{\longrightarrow}\pi_1 (V_{n+1,n}) \to 0\qquad,$$
obtained from the long homotopy exact sequence. After all, ``everyone knows'' that $\SS^n$ is classically aspherical when $n$ is greater than 2.

There are conceptual differences in the calculation of the fundamental-group $\pi_1(\SS^n)$ compared with that of $\pi_2 (\SS^n)$. A simplified version of the Seifert-van Kampen Theorem given in  \cite{Munkres} shows that $\SS^n$ is simply connected. Even simpler, Proposition 1.14 of \cite{Hatcher} shows how to homotope any loop $\gamma$ in $\SS^n$ to avoid a given value $x_0 \in \SS^n$. Thus the new loop actually resides in a space topologically equivalent to the contractible $\RR^n$, hence shrinks to a constant loop.

The modifications we saw  in the loop $\gamma$ were local perturbations. A standard proof starting with
$\sigma : \left(\SS^2, x_0\right) \to \left(\SS^n, y_0\right)$ wants to perform the same kind of  approximation on a representative $\sigma$ of element $[\sigma] \in \pi_2 \left(\SS^n, y_0\right)$. A close linear approximation $\hat{\sigma}$ of $\sigma$ will work since the two neighboring ``immersed'' 2-spheres can never be antipodal on $\SS^n$, hence must be 2-loops homotopic to each other, see \cite{Dugundji}. The approximation $\hat{\sigma}$, uncovering at least one image point, is homotopic to $\sigma$ by a {\it global} argument. Besides this, the ``linear approximation'' construction features an implicit notion of ``oriented simplex''. But we are actually researching the under-pinnings of ``oriented simplex'', ``determinant'', and ``parity of an $n$-frame'', so it is more consistent to avoid these ``orientation'' concepts from our proofs.

In \cite{Althoen}, the author develops a ``van Kampen Theorem'' for $\pi_2$. In the case of $\SS^n$, one decomposes the space into two open hemi-spheres, which overlap in a (simply connected) thickened equator. The sought-after result that
$\pi_2 (\SS^n) = 0$ follows from a Hurewicz theorem which is almost a tautology if one employs {\it cubical} homology, see \cite{Hatcher}. But for this, one must lay the foundations of an entire homological theory, which already brings in the concept of ``oriented cube''.

For us, none of these simplicial or homological approaches are quite satisfactory, and we are left with methods that are more algebraic in nature \cite{Althoen}, \cite{Brown}.

Althoen shows that $\pi_2$ obeys a commutative push-out diagram of the form
$$\left.\matrix
\pi_2(Y) & \rightarrow & \pi_2 (A) \\
\downarrow & & \downarrow \\
\pi_2(B) & \rightarrow & \pi_2 (\SS^n)
\endmatrix\right.$$
where $Y = A \cap B$, and all mappings are induced from inclusions. Such a ``push-out'' illustrates how $\pi_2 (\SS^n)$ is generated by $\pi_2(A)$ and $\pi_2 (B)$, but both of these homotopy groups are trivial, hence $\pi_2 (\SS^n) = 0$, the zero group. We observed that Althoen's method relies on the Hurewicz Isomorphism Theorem, or ``Blakers-Massey'' result \cite{Hatcher}, so we seek another approach that does not bring in the algebra of chains, cycles and boundaries, i.e., the homology (with integer coefficients) of cell complexes.

The approach toward $\pi_2(\SS^n) = 0$ that we prefer, comes from the theory of ``crossed modules'' originally developed by \cite{Brown, Higgins}. We gain results analogous to a ``van Kampen'' Theorem that applies to $\pi_2$, with the flavor of ``geometric group theory''. For example, it is proven in this theory that for $Z$ a path-connected space, the 1-fold {\it suspension} $S Z$ satisfies
$$\pi_2 (S Z) = \pi_1 (Z)\;\;\;\text{``abelianized''}\qquad.$$
In our cases of interest $Z = \SS^2\!, \SS^3\!, \dotsc, $ where we know $\pi_1(Z) = 0$, we may conclude that $\pi_2 (\SS^3) = 0$, $\pi_2(\SS^4) = 0$ and so forth. Thus the Brown-Higgins result on $\pi_2$ of a {\it suspension} seems appropriate to our study of the parity of 3-frames, 4 frames, ...~.

In order to determine the path-components of $O(n) \simeq V_{n,n}$ for $n \geq 3$, the critical calculation was to find $\pi_1 (V_{n,n-1})$. This comes out to be the group of two elements, based on the ``special case'' of homotopy, namely that $\pi_1 (\SS^{n-1})= 0$. If we take the position that we {\it do not} have a handle on $\pi_2(\SS^m)$ for large $m$, we may only be permitted to conclude that $\pi_1 (V_{j,j-1}) \simeq \ZZ_2$, $j \leq m$, but that $\pi_1 (V_{m+1,m}) = 0$, and so for all indices $m' \geq  m$. In another sense, the calculation of $\pi_2$ is a {\it special case} as well. For $X = \SS^{n-1}$ ``simply connected'' and homologically trivial $(n > 3)$, the Hurewicz isomorphism
$$h: \pi_2 (X) \to H_2(X)$$
gives the result $\pi_2(X) = 0$ which is needed for the ``strong'' calculation of $\pi_1(V_{n,n-1})$\linebreak
$ = \ZZ_2$. Using say ``cubical'' homology, the mapping $h$ is straightforward to define: for $f: I^2 \to X$ with $f(\partial I^2) = x_0$, base point, let $h\left([f]\right) = \{f\}$, the latter being the cycle class of the singular 2-cube $[f]$, modulo integral boundary chains. Since $H_2(X) = 0$ for spaces within our present interest, not a lot of homological calculation will be necessary. But we do need to consider equivalences in $H_2(X; \ZZ)$, which means that {\it cubical} 3-chains on $X$ should not be ignored. In any case it is possible to define a left inverse $r: H_2(X; \ZZ) \to \pi_2 (X)$ in general, so that $r \circ h(s) = s$, indicating that $h$ is an injection for $X = \SS^j$, $j \geq 2$, and $\pi_2(X) = 0$ for $j > 2$.

We accepted, provisionally, the Hurewicz ``isomorphism theorem'' under a Principle of Low Dimension. This underlies Althoen's version of a van Kampen Theorem for $\pi_2$. Also, using a mapping $g : \SS^2 \to X$ as the {\it initial data} for a heat equation on $X$: the ``heat flow'' instantly fabricates a homotopy $g_t$ from $g$ to a smooth mapping $g_{t_1} : \SS^2 \to X$, where $t_1 > 0$. The latter mapping must miss some value in $X = \SS^{n-1}$ according to the theorem of A\. B\. Brown, or that of Sard, \cite{Milnor, TFDV}, and hence be contractible to a constant (point) mapping.

All in all, it seems justified to adopt {\it  one} of these arguments for the ``strong case'' that, actually, ·$\pi_1 (V_{n,n-1}) \simeq \ZZ$, hence reaching the conclusion
$$\pi_1(V_{n,n-1}) \simeq \pi_1 (V_{n+1,n}) \simeq \ZZ_2$$
for $n \geq 3$. We summarize in the following Section the salient point that two distinct elements (or $n$-frame ``parity'') occur in $O(n)$, namely ``rotations'' and rotations-with-reflection.

\head The Compact Quaternion Cover\endhead

The only ``universal'' covering that need be mentioned is a double covering of $V_{3,2} \simeq O^*(3)$. The isomorphism on $\pi_1$ that is derived from ``strong form'' of the exact homotopy sequence yields
$$\pi_1 (V_{3,2}) \simeq  \pi_1 (V_{4,3}) \simeq \cdots \simeq \pi_1 (V_{n+1,n}) \simeq \ZZ_2\qquad.$$
The long homotopy sequence also provides canonical generating elements $\xi_k \in \pi_1 (V_{k, k-1})$. In fact, a loop called the canonical Givens loop $\gamma_k$ represents $\xi_k$, and is easily visualized by means of a non-contractible loop on $V_{3,2}$.

The real algebra of quaternions is a $4$-dimensional vector space $\HH$, often seen as generated by {\bf 1},  $\pmb{i}$, $\pmb{j}$, $\pmb{k}$. The ``3-sphere'' $\SS^3$ or group of unit quaternions $\HH_1$, is given by coordinates $(x, y, z, w)$ such that
$$
x^2 + y^2 + z^2 + w^2 = 1\qquad. \tag{*}
$$
Bearing in mind the condition (*), a continuous and algebraic mapping
$$\Phi: \HH_1 \to O(3)\qquad,$$
not necessarily surjective, can be given by
$$\Phi(x,y,z,w) = \left.\matrix
x^2+y^2-z^2-w^2 & 2yz-2xw & 2yw+2xz \\
2yz+2xw & x^2-y^2+z^2-w^2 & 2zw-2xy \\
2yw-2xz & 2zw+2xy & x^2-y^2-z^2+w^2
\endmatrix\right. \qquad. $$
If $q = (x, y, z, w)$ and $\sigma = (y, z, w) \in \RR^3$, then {\it Rodrigues' formula} is derivable from
$\kappa = (\alpha, \beta,\gamma) \in\RR^3$
$$\Phi (x,y,z,w) \cdot (\alpha, \beta,\gamma) = \left(\alpha', \beta', \gamma'\right) \tag{1.1}$$
where $\kappa' = \kappa +2x (\sigma \times \kappa) + 2(\sigma \times (\sigma \times \kappa))$, recalling that $x$ is ``the {\it scalar part}'' of quaternion $q$.

We also keep in mind the definition of (unit) quaternion multiplication as derived from
$$q'' = q \cdot q'$$
or
$$\align
x'' &= xx' - yy' -zz' -ww' \\
y'' & = xy' + yx' -zw'-wz' \\
z'' &= xz' +zx' -wy'+yw' \\
w'' &= xw' +wx' -yz'+zy' \qquad.
\endalign$$
That the quaternion product should be well-defined requires that
$$\|q''\|^2 = \left(x''\right)^2 + \left(y''\right)^2 + \left(z''\right)^2 + \left(w''\right)^2\qquad.$$

We use automated symbolic algebra from the Unix Octave package to verify the important properties of $\HH_1$ and the mapping $\Phi$. One must check that $v, w\in \HH$, of norm $1$, give a product quaternion $vw$ that is also a unit 4-vector (has ``norm squared'' $= 1)$.

This fact is shown in the edited QuarterDoc (as part of a Figure).

Next, it would be nice to see that the {\it product} of $\HH_1$ is associative. But this is nearly obvious in view of the construction of the vector space $\HH$ and the generating rules of multiplication $i \cdot i = -1$, $j \cdot j = -1$, $k \cdot -k = -1$, $i \cdot j = k$ and so forth. However, one may demonstrate the associativity on $\HH_1$ directly (Aa\_verify). Then the ``inverse law'' should be checked, namely
$$q^{-1} = q^*/\|q^*\|^2$$
where $q^*$ is the conjugate $(q_1, -q_2, -q_3, -q_4)$.

But it is also obvious that the conjugate of a unit norm quaternion also has unit norm, so in fact $q \in \HH_1$ implies that $q^{-1} = q^*$.

Among necessary properties of the mapping $\Phi$ are that it is well-defined. In other words when $q$ has unit norm, $\phi(q) \in O(3)$. That is, we get a $3 \times 3$ matrix with pairwise orthogonal rows of unit norm.

This is demonstrated in the file  Aa\_verify.

We take it as known that $\SS^3 \simeq \HH_1$ as well as $O(3)$, is a smooth compact manifold, in fact both are Lie groups. The mapping $\Phi: \HH_1 \to O(3)$ is also (infinitely) smooth, as defined by polynomial entries. The key topological result, that $\Phi$ is a finite covering (in fact two-to-one) onto its image, follows transparently if we are allowed to use the very concrete group theory that lies at hand.

So it is necessary to show that $\Phi$ is a homomorphism of groups. It is easy to see that $\pmb{1}$ goes to the identity
$$\Phi(\pmb{1}) = \left[\matrix
1 & 0 & 0 \\
0 & 1 & 0 \\
0 & 0 & 1
\endmatrix\right]$$
the {\it trivial rotation} of three-space.
\vskip-12pt

\newpage

{{\tt



\%\%  Show that the quaternion product of two Norm One quaternions has Norm One

\bigskip

>> syms v [1 4] real

>> norm(v)

 ans =

 (abs(v1)\^{}2 + abs(v2)\^{}2 + abs(v3)\^{}2 + abs(v4)\^{}2)\^{}(1/2)

>> nrmsq(v)

ans =

v1\^{}2 + v2\^{}2 + v3\^{}2 + v4\^{}2

\bigskip

>> [f g h k]=qprod(v,w)

f =

v1*w1 - v2*w2 - v3*w3 - v4*w4

g =

v1*w2 + v2*w1 + v3*w4 - v4*w3

h =

v1*w3 + v3*w1 - v2*w4 + v4*w2

k =

v1*w4 + v2*w3 + v3*w2 + v4*w1

>> zz= [f g h k]

zz =

[ v1*w1 - v2*w2 - v3*w3 - v4*w4, v1*w2 + v2*w1 + v3*w4 - v4*w3,

v1*w3 + v3*w1 - v2*w4 + v4*w2, v1*w4 + v2*w3 + v3*w2 + v4*w1]

>> nfz=nrmsq(zz)

nfz =

(v1*w4 + v2*w3 + v3*w2 + v4*w1)\^{}2 + (v1*w2 + v2*w1 + v3*w4 - v4*w3)\^{}2 +

(v1*w3 + v3*w1 - v2*w4 + v4*w2)\^{}2 + (v2*w2 - v1*w1 + v3*w3 + v4*w4)\^{}2

>> assume(nrmsq(v)==1)

>> assume(nrmsq(w)==1)

>> simplify(nfz)

ans =

1

}

\medskip
\centerline{Norm Product}


\newpage

{{\tt

\bigskip

Homomorphism Phi(a) = Aa (Aa\_{}verify)

The definition of the homomorphism starting with the quaternion    a     gives a 3  by   3 real matrix

>> Aa=so3(a)

Aa =

[ a1\^{}2 + a2\^{}2 - a3\^{}2 - a4\^{}2,         2*a2*a3 - 2*a1*a4,         2*a1*a3 + 2*a2*a4]

[         2*a1*a4 + 2*a2*a3, a1\^{}2 - a2\^{}2 + a3\^{}2 - a4\^{}2,         2*a3*a4 - 2*a1*a2]

[         2*a2*a4 - 2*a1*a3,         2*a1*a2 + 2*a3*a4, a1\^{}2 - a2\^{}2 - a3\^{}2 + a4\^{}2]

Now we see tha Aa '  *  Aa is the identity

>> assume(a,'real')

>> Asq=Aa*Aa'

Asq =

[                                                               (2*a1*a3 + 2*a2*a4)\^{}2 + (2*a1*a4 - 2*a2*a3)\^{}2 + (a1\^{}2 + a2\^{}2 - a3\^{}2 - a4\^{}2)\^{}2, (2*a1*a4 + 2*a2*a3)*(a1\^{}2 + a2\^{}2 - a3\^{}2 - a4\^{}2) - (2*a1*a4 - 2*a2*a3)*(a1\^{}2 - a2\^{}2 + a3\^{}2 - a4\^{}2) - (2*a1*a2 - 2*a3*a4)*(2*a1*a3 + 2*a2*a4), (2*a1*a3 + 2*a2*a4)*(a1\^{}2 - a2\^{}2 - a3\^{}2 + a4\^{}2) - (2*a1*a3 - 2*a2*a4)*(a1\^{}2 + a2\^{}2 - a3\^{}2 - a4\^{}2) - (2*a1*a2 + 2*a3*a4)*(2*a1*a4 - 2*a2*a3)]

[ (2*a1*a4 + 2*a2*a3)*(a1\^{}2 + a2\^{}2 - a3\^{}2 - a4\^{}2) - (2*a1*a4 - 2*a2*a3)*(a1\^{}2 - a2\^{}2 + a3\^{}2 - a4\^{}2) - (2*a1*a2 - 2*a3*a4)*(2*a1*a3 + 2*a2*a4),                                                               (2*a1*a2 - 2*a3*a4)\^{}2 + (2*a1*a4 + 2*a2*a3)\^{}2 + (a1\^{}2 - a2\^{}2 + a3\^{}2 - a4\^{}2)\^{}2, (2*a1*a2 + 2*a3*a4)*(a1\^{}2 - a2\^{}2 + a3\^{}2 - a4\^{}2) - (2*a1*a2 - 2*a3*a4)*(a1\^{}2 - a2\^{}2 - a3\^{}2 + a4\^{}2) - (2*a1*a3 - 2*a2*a4)*(2*a1*a4 + 2*a2*a3)]

[ (2*a1*a3 + 2*a2*a4)*(a1\^{}2 - a2\^{}2 - a3\^{}2 + a4\^{}2) - (2*a1*a3 - 2*a2*a4)*(a1\^{}2 + a2\^{}2 - a3\^{}2 - a4\^{}2) - (2*a1*a2 + 2*a3*a4)*(2*a1*a4 - 2*a2*a3), (2*a1*a2 + 2*a3*a4)*(a1\^{}2 - a2\^{}2 + a3\^{}2 - a4\^{}2) - (2*a1*a2 - 2*a3*a4)*(a1\^{}2 - a2\^{}2 - a3\^{}2 + a4\^{}2) - (2*a1*a3 - 2*a2*a4)*(2*a1*a4 + 2*a2*a3),                                                               (2*a1*a2 + 2*a3*a4)\^{}2 + (2*a1*a3 - 2*a2*a4)\^{}2 + (a1\^{}2 - a2\^{}2 - a3\^{}2 + a4\^{}2)\^{}2]

>> size(Asq)

     3     3

>> simplify(expand(Asq))

ans =

[ (a1\^{}2 + a2\^{}2 + a3\^{}2 + a4\^{}2)\^{}2,                             0,                             0]

[                             0, (a1\^{}2 + a2\^{}2 + a3\^{}2 + a4\^{}2)\^{}2,                             0]

[                             0,                             0, (a1\^{}2 + a2\^{}2 + a3\^{}2 + a4\^{}2)\^{}2]

>> Aans=ans;

>> assume(nrmsq(a)==1)  

>> simplify(Aans)

ans =

[ 1, 0, 0]

[ 0, 1, 0]

[ 0, 0, 1]

Hence     Aa represents    a rotation in   3-space

}

\medskip
\centerline{Image is an Orthogonal Motion}


\newpage

{{\tt

\%\% Calculations to show that the quaternion product is Associative
\bigskip

**  Short Version of Associativity
\bigskip

v = qprod(b,c)

[ b1*c1 - b2*c2 - b3*c3 - b4*c4, b1*c2 + b2*c1 + b3*c4 - b4*c3,

b1*c3 + b3*c1 - b2*c4 + b4*c2, b1*c4 + b2*c3 - b3*c2 + b4*c1]

 >> [w1 w2 w3 w4]=qprod(a,v)

w1 =

- a1*(b2*c2 - b1*c1 + b3*c3 + b4*c4) - a2*(b1*c2 + b2*c1 + b3*c4 - b4*c3) -

a3*(b1*c3 + b3*c1 - b2*c4 + b4*c2) - a4*(b1*c4 + b2*c3 - b3*c2 + b4*c1)

>> [u1 u2 u3 u4] = qprod(a,b)

u1 =

a1*b1 - a2*b2 - a3*b3 - a4*b4

>> z = [z1 z2 z3 z4]=qprod(u,c)

z1=

- c1*(a2*b2 - a1*b1 + a3*b3 + a4*b4) - c2*(a1*b2 + a2*b1 + a3*b4 - a4*b3)

- c3*(a1*b3 + a3*b1 - a2*b4 + a4*b2) - c4*(a1*b4 + a2*b3 - a3*b2 + a4*b1)

>> z1-w1 = ans

a1*(b2*c2 - b1*c1 + b3*c3 + b4*c4) - c2*(a1*b2 + a2*b1 + a3*b4 - a4*b3)

 - c3*(a1*b3 + a3*b1 - a2*b4 + a4*b2) - c4*(a1*b4 + a2*b3 - a3*b2 + a4*b1)

- c1*(a2*b2 - a1*b1 + a3*b3 + a4*b4) + a2*(b1*c2 + b2*c1 + b3*c4 - b4*c3)

+ a3*(b1*c3 + b3*c1 - b2*c4 + b4*c2) + a4*(b1*c4 + b2*c3 - b3*c2 + b4*c1)

>> simplify(ans)

 ans =   0

}

\medskip
\centerline{Quaternion Product Associativity Law}


\newpage

{{\tt

\bigskip

Homomorphism Property of the Quaternion Mapping  Phi  to the Rotation Group (Homomorph)

>> Aa=so3(a)

[ a1\^{}2 + a2\^{}2 - a3\^{}2 - a4\^{}2,         2*a2*a3 - 2*a1*a4,         2*a1*a3 + 2*a2*a4]

[         2*a1*a4 + 2*a2*a3, a1\^{}2 - a2\^{}2 + a3\^{}2 - a4\^{}2,         2*a3*a4 - 2*a1*a2]

[         2*a2*a4 - 2*a1*a3,         2*a1*a2 + 2*a3*a4, a1\^{}2 - a2\^{}2 - a3\^{}2 + a4\^{}2]

 >> Ab=so3(b)

 [ b1\^{}2 + b2\^{}2 - b3\^{}2 - b4\^{}2,         2*b2*b3 - 2*b1*b4,         2*b1*b3 + 2*b2*b4]

[         2*b1*b4 + 2*b2*b3, b1\^{}2 - b2\^{}2 + b3\^{}2 - b4\^{}2,         2*b3*b4 - 2*b1*b2]

[         2*b2*b4 - 2*b1*b3,         2*b1*b2 + 2*b3*b4, b1\^{}2 - b2\^{}2 - b3\^{}2 + b4\^{}2]

>> LAprod=Aa*Ab;

>> Ac=so3(c)

 [ (a1*b2 + a2*b1 + a3*b4 - a4*b3)\^{}2 - (a1*b3 + a3*b1 - a2*b4 + a4*b2)\^{}2 - (a1*b4 + a2*b3 - a3*b2 + a4*b1)\^{}2 + (a2*b2 - a1*b1 + a3*b3 + a4*b4)\^{}2,         2*(a1*b2 + a2*b1 + a3*b4 - a4*b3)*(a1*b3 + a3*b1 - a2*b4 + a4*b2) + 2*(a1*b4 + a2*b3 - a3*b2 + a4*b1)*(a2*b2 - a1*b1 + a3*b3 + a4*b4),         2*(a1*b2 + a2*b1 + a3*b4 - a4*b3)*(a1*b4 + a2*b3 - a3*b2 + a4*b1) - 2*(a1*b3 + a3*b1 - a2*b4 + a4*b2)*(a2*b2 - a1*b1 + a3*b3 + a4*b4)]

[         2*(a1*b2 + a2*b1 + a3*b4 - a4*b3)*(a1*b3 + a3*b1 - a2*b4 + a4*b2) - 2*(a1*b4 + a2*b3 - a3*b2 + a4*b1)*(a2*b2 - a1*b1 + a3*b3 + a4*b4), (a1*b3 + a3*b1 - a2*b4 + a4*b2)\^{}2 - (a1*b2 + a2*b1 + a3*b4 - a4*b3)\^{}2 - (a1*b4 + a2*b3 - a3*b2 + a4*b1)\^{}2 + (a2*b2 - a1*b1 + a3*b3 + a4*b4)\^{}2,         2*(a1*b2 + a2*b1 + a3*b4 - a4*b3)*(a2*b2 - a1*b1 + a3*b3 + a4*b4) + 2*(a1*b3 + a3*b1 - a2*b4 + a4*b2)*(a1*b4 + a2*b3 - a3*b2 + a4*b1)]

[         2*(a1*b2 + a2*b1 + a3*b4 - a4*b3)*(a1*b4 + a2*b3 - a3*b2 + a4*b1) + 2*(a1*b3 + a3*b1 - a2*b4 + a4*b2)*(a2*b2 - a1*b1 + a3*b3 + a4*b4),         2*(a1*b3 + a3*b1 - a2*b4 + a4*b2)*(a1*b4 + a2*b3 - a3*b2 + a4*b1) - 2*(a1*b2 + a2*b1 + a3*b4 - a4*b3)*(a2*b2 - a1*b1 + a3*b3 + a4*b4), (a1*b4 + a2*b3 - a3*b2 + a4*b1)\^{}2 - (a1*b3 + a3*b1 - a2*b4 + a4*b2)\^{}2 - (a1*b2 + a2*b1 + a3*b4 - a4*b3)\^{}2 + (a2*b2 - a1*b1 + a3*b3 + a4*b4)\^{}2]

>> NonHomom=LAprod-Ac

 [ (a1*b3 + a3*b1 - a2*b4 + a4*b2)\^{}2 - (a1*b2 + a2*b1 + a3*b4 - a4*b3)\^{}2 + (a1*b4 + a2*b3 - a3*b2 + a4*b1)\^{}2 - (a2*b2 - a1*b1 + a3*b3 + a4*b4)\^{}2 + (a1\^{}2 + a2\^{}2 - a3\^{}2 - a4\^{}2)*(b1\^{}2 + b2\^{}2 - b3\^{}2 - b4\^{}2) - (2*a1*a3 + 2*a2*a4)*(2*b1*b3 - 2*b2*b4) - (2*a1*a4 - 2*a2*a3)*(2*b1*b4 + 2*b2*b3),         (2*a1*a3 + 2*a2*a4)*(2*b1*b2 + 2*b3*b4) - (2*b1*b4 - 2*b2*b3)*(a1\^{}2 + a2\^{}2 - a3\^{}2 - a4\^{}2) - 2*(a1*b2 + a2*b1 + a3*b4 - a4*b3)*(a1*b3 + a3*b1 - a2*b4 + a4*b2) - 2*(a1*b4 + a2*b3 - a3*b2 + a4*b1)*(a2*b2 - a1*b1 + a3*b3 + a4*b4) - (2*a1*a4 - 2*a2*a3)*(b1\^{}2 - b2\^{}2 + b3\^{}2 - b4\^{}2),         (2*a1*a3 + 2*a2*a4)*(b1\^{}2 - b2\^{}2 - b3\^{}2 + b4\^{}2) + (2*b1*b3 + 2*b2*b4)*(a1\^{}2 + a2\^{}2 - a3\^{}2 - a4\^{}2) - 2*(a1*b2 + a2*b1 + a3*b4 - a4*b3)*(a1*b4 + a2*b3 - a3*b2 + a4*b1) + 2*(a1*b3 + a3*b1 - a2*b4 + a4*b2)*(a2*b2 - a1*b1 + a3*b3 + a4*b4) + (2*a1*a4 - 2*a2*a3)*(2*b1*b2 - 2*b3*b4)]

[         (2*a1*a4 + 2*a2*a3)*(b1\^{}2 + b2\^{}2 - b3\^{}2 - b4\^{}2) + (2*b1*b4 + 2*b2*b3)*(a1\^{}2 - a2\^{}2 + a3\^{}2 - a4\^{}2) - 2*(a1*b2 + a2*b1 + a3*b4 - a4*b3)*(a1*b3 + a3*b1 - a2*b4 + a4*b2) + 2*(a1*b4 + a2*b3 - a3*b2 + a4*b1)*(a2*b2 - a1*b1 + a3*b3 + a4*b4) + (2*a1*a2 - 2*a3*a4)*(2*b1*b3 - 2*b2*b4), (a1*b2 + a2*b1 + a3*b4 - a4*b3)\^{}2 - (a1*b3 + a3*b1 - a2*b4 + a4*b2)\^{}2 + (a1*b4 + a2*b3 - a3*b2 + a4*b1)\^{}2 - (a2*b2 - a1*b1 + a3*b3 + a4*b4)\^{}2 + (a1\^{}2 - a2\^{}2 + a3\^{}2 - a4\^{}2)*(b1\^{}2 - b2\^{}2 + b3\^{}2 - b4\^{}2) - (2*a1*a2 - 2*a3*a4)*(2*b1*b2 + 2*b3*b4) - (2*a1*a4 + 2*a2*a3)*(2*b1*b4 - 2*b2*b3),         (2*a1*a4 + 2*a2*a3)*(2*b1*b3 + 2*b2*b4) - (2*b1*b2 - 2*b3*b4)*(a1\^{}2 - a2\^{}2 + a3\^{}2 - a4\^{}2) - 2*(a1*b2 + a2*b1 + a3*b4 - a4*b3)*(a2*b2 - a1*b1 + a3*b3 + a4*b4) - 2*(a1*b3 + a3*b1 - a2*b4 + a4*b2)*(a1*b4 + a2*b3 - a3*b2 + a4*b1) - (2*a1*a2 - 2*a3*a4)*(b1\^{}2 - b2\^{}2 - b3\^{}2 + b4\^{}2)]

[         (2*a1*a2 + 2*a3*a4)*(2*b1*b4 + 2*b2*b3) - (2*b1*b3 - 2*b2*b4)*(a1\^{}2 - a2\^{}2 - a3\^{}2 + a4\^{}2) - 2*(a1*b2 + a2*b1 + a3*b4 - a4*b3)*(a1*b4 + a2*b3 - a3*b2 + a4*b1) - 2*(a1*b3 + a3*b1 - a2*b4 + a4*b2)*(a2*b2 - a1*b1 + a3*b3 + a4*b4) - (2*a1*a3 - 2*a2*a4)*(b1\^{}2 + b2\^{}2 - b3\^{}2 - b4\^{}2),         (2*a1*a2 + 2*a3*a4)*(b1\^{}2 - b2\^{}2 + b3\^{}2 - b4\^{}2) + (2*b1*b2 + 2*b3*b4)*(a1\^{}2 - a2\^{}2 - a3\^{}2 + a4\^{}2) + 2*(a1*b2 + a2*b1 + a3*b4 - a4*b3)*(a2*b2 - a1*b1 + a3*b3 + a4*b4) - 2*(a1*b3 + a3*b1 - a2*b4 + a4*b2)*(a1*b4 + a2*b3 - a3*b2 + a4*b1) + (2*a1*a3 - 2*a2*a4)*(2*b1*b4 - 2*b2*b3), (a1*b2 + a2*b1 + a3*b4 - a4*b3)\^{}2 + (a1*b3 + a3*b1 - a2*b4 + a4*b2)\^{}2 - (a1*b4 + a2*b3 - a3*b2 + a4*b1)\^{}2 - (a2*b2 - a1*b1 + a3*b3 + a4*b4)\^{}2 + (a1\^{}2 - a2\^{}2 - a3\^{}2 + a4\^{}2)*(b1\^{}2 - b2\^{}2 - b3\^{}2 + b4\^{}2) - (2*a1*a2 + 2*a3*a4)*(2*b1*b2 - 2*b3*b4) - (2*a1*a3 - 2*a2*a4)*(2*b1*b3 + 2*b2*b4)]

 >> simplify(expand(NonHomom))

 [ 0, 0, 0]

[ 0, 0, 0]

[ 0, 0, 0]

}

\medskip
\centerline{Three Dimensional Homomorphism}


\newpage

Then we need the property
$$\Phi(q \cdot q') = \Phi (q) \Phi (q'), \tag{3.1}$$
where group multiplication on the Right-Hand side is effectively {\it linear composition} (matrix multiplication). It is seen in Homomorph that the entries of the matrices (3.1) are the same. For an exemplar we look at LAprod(1,1)-Ac(1,1) in this file. This gives the (1,1) entry of the matrix, the difference between the {\it matrix product} of $Aa = \Phi(a)$ and $Ab = \Phi(b)$, and the orthogonal matrix $Ac = \Phi(c)$ corresponding to the quaternion product $a \cdot b$. This entry, as well as the other entries of this $3 \times 3$ difference matrix, is $0$.

It is of interest to examine the ``fiber'' of $\Phi$ at $I_{3 \times 3} \in O(3)$. This is seen ``by hand'' to consist of the two ``unit'' quaternions $e = (1, 0, 0, 0)$ and $-e = (-1,0,0,0)$. Of course $e$ is also the unity element of the group $\HH_1$. Now we apply the Principle of Low Dimension in the guise of our known Lie groups, and pose $\Phi(a) = \Phi(b)$. But by the properties of a ``homeomorphism'', we may derive $\Phi(a b^{-1}) = I_{3\times 3}$, so $a \cdot b^{-1} = e$ or $-e$. Hence, as ``unit'' quaternions, $a = b$ or $a = -b$. Hence $\Phi$ is a true double covering, without ``branch points''.

To look ahead, it is unsurprising that our major effort in analyzing frames of $m$-vectors occurs in dimension three. From our simple reasoning thus far, we obtain a result arguably of genuine mathematics.

Notice that the 3-frame
$$\pmb{f} = \left[\matrix
1 & 0 & 0 \\
0 & 1 & 0 \\
0 & 0 & -1
\endmatrix\right]$$
cannot be in the image of $\Phi$. But $\pmb{f} \in O(3)$. According to Brouwer's ``Invariance of Domain'' factum, a ``locally finite'' covering $G: N \to M$, where $N$ and $M$ are compact 3-manifolds (or $k$-manifolds) with $M$ {\it path-connected}, must be a {\it surjection}, \cite{Hatcher}, p. 172. See also \cite{Sjogren, Homogeneous} and \cite{Sjogren, Domain}. Therefore $O(3)$ is {\it not} path-connected, so the 3-frame $\pmb{f}$ has no ``orthogonal path'' to the identity frame that we designated as $I_{3 \times 3}$.

A specific construction will be given the next Section, showing that $\Phi$ on $\HH_1 \simeq \SS^3$ maps {\it onto} $V_{3,2}$ which we see must form a path-component of $O(3)$, naturally {\it omitting} the 3-frame called $\pmb{f}$.

We have solved the parity problem for the topological space $O(3) \simeq V_{3 \cdot 3}$, showing that there are exactly two path-components. This conclusion will also follow if we examine the action, by deck transformation, of the double covering $\rho_3: O(3) \to V_{3,2}$. Since for $n > 3$, the study of the universal covering
$$\Lambda _n : \text{Spin}(n) \to V_{n, n-1}$$
is more intricate than that of our group homomorphism $\Phi: \SS^3 \to O^+(3)$, we use the homotopy exact sequence of the canonical fibration to calculate $\pi_1 (V_{n, n-1})$, as indicated above.

The remaining point of importance is to construct a geometric generating loop for $\pi_1$. We may consider 4-vectors in $\HH_1$ as follows.
$$\align
q_0 &= (1,0,0,0)  \quad {\Cal N} \quad \text{(North Pole)} \\
q_{\theta} &= (\cos \theta, \sin \theta, 0, 0)  \\
q_{\pi} &= (-1, 0,0,0)  \quad {\Cal S} \quad \text{(South Pole)}
\endalign$$
so that
$$\Phi q_0 = \left[\matrix
1 & \cdots & \\
&  1 & \cdots \\
&& 1
\endmatrix\right]$$
and
$$\Phi q_{\theta} = \left[\matrix
1 & 0 & 0 \\
0 & \cos 2\theta & -\sin 2\theta \\
0 &\sin 2\theta & \cos 2\theta
\endmatrix\right],\qquad \Phi q_{\pi} = \left[\matrix
1 & 0 & 0 \\
 & 1 & 0 \\
 && 1
 \endmatrix\right]\quad \text{also.}$$
 It seems clear that $\gamma(\theta) = \Phi q_{\theta}$ is some kind of {\it Givens loop} on $V_{3,2}$. This loop must be a generator of $\pi_1 (O^+(3)) \simeq \ZZ_2$. It cannot be contractible to a point since the loop is covered by the {\it path} $q_{\theta}$, $0 \leq \theta \leq \pi$, which has distinct end-points ${\Cal N}$ and ${\Cal S}$.

 We may return to consideration of the ``strong case'' of a segment of the long homotopy exact sequence
$$\align
  \longrightarrow &\pi_2 (\SS^n)  \buildrel\hbox{$\partial$}\over{\longrightarrow}\pi_1 (V_{n,n-1})  \buildrel\hbox{$\iota_{*}$}\over{\longrightarrow} \pi_1 (V_{n+1,n}) \\
 & \buildrel\hbox{$\rho^{n+1}_{*}$}\over{\longrightarrow} \pi_1 (\SS^n) = 0\qquad,
\endalign$$
where $n \geq 3$ and $\rho^{n+1}: V_{n+1,n} \to \SS^n$ gives the {\it normalized} vector in $\RR^{n+1}$ which amounts to the {\it first} vector $w_1$ of the given $n$-frame $\pmb{w} = \left(w_1, \dotsc, w_n\right)$. The ``strong condition'' refers to taking $\pi_2 (\SS^n) = 0$ via one of various proof methodologies and arrive at
$$\pi_1 \left(V_{n, n-1}\right) \buildrel\hbox{$\iota_{*}$}\over{\longrightarrow} \pi_1 \left(V_{n+1,n}\right)$$
as an isomorphism. We already know that $\pi_1 (V_{3,2}) \simeq \ZZ_2$, generated by $\gamma (\theta) = \Phi q_{\theta}$, a Givens loop on the $y,z$ plane consisting of a continuous rotation from $\phi=0$ to $\phi = 2\pi$ that keeps {\it fixed} the $x$-axis.

The {\it proof} of the long homotopy exact sequence of a linear bundle expresses the fact that $\pi_1 (F) \to \pi_1 (E)$ is induced by the inclusion of the {\it fiber} $F$ into the total space $E$. Thus the generator $\gamma (\theta)$ of $\pi_1 (V_{3,2})$ must go by $\iota_*$ to a generator of $\pi_1 (V_{4,3})$.

Related material is covered in Theorem 4.41 of \cite{Hatcher}. The latter generator is provided by the ``same'' Givens loop, or rather a suspension of $\gamma (\theta)$, one which now leaves fixed the $x$ and $w$ axes of $\RR^4$. We have at this point completely set up an induction step for $\pi_1 (O^+(n))$, giving a geometric generator of the same.

\bigskip
\noindent
{\bf Proposition D}\quad  The fundamental group of the compact space of $(n-1)$-frames in $\RR^n$ equals the two-group $\ZZ_2$. It has as generator $\gamma^s(\theta)$, $0 \le \theta \le \pi$, the continuous and uniform rotation of the $y,z$ plane, suspended sufficiently many times to a continuous path of $O(n)$ as
$$\gamma^s(\theta) = \left.\matrix
1 & 0 & 0 & 0 & \cdots & 0 \\
0 & \cos 2\theta &-\sin 2\theta & \\
0 & -\sin 2\theta & \cos 2 \theta & \\
\vdots & \hfil&&\ddots &  \\
& & && 1 & \\
0 & \hfil & \cdots && & 1 &.
\endmatrix\right. \tag{G}$$
Finally, we may set down the result that two $n$-frames of opposite parity such as $\pmb{e} = (e_1, e_2, e_3, \dotsc, e_n)$ and $\pmb{e'} = (e_2, e_1, e_3,\dotsc, e_n)$ cannot be connected through a continuous path in $O(n) = V_{n,n}$. Our context applies to the compact orthogonal ``group'', on which the ``determinant'' assumes only a discrete set of real values.

\bigskip
\noindent
{\bf Theorem}\quad For $n \geq 3$, the double covering $\rho_n = O(n) \simeq V_{n,n} \to V_{n,n-1}$ maps a non-path-connected space onto a path-connected space (maximal component).

\bigskip
\noindent
{\it Proof.} \quad We already saw that $\rho_n$ is surjective (the deletion of the final vector $v_n$ from an $n$-frame) and that $V_{n,n-1}$ is path-connected. Suppose that $O(n)$ were path-connected. Then there would be a path $\delta (t)$ in $O(n)$ from $\pmb{e}$ to $\pmb{\overline{e}} = (e_1, e_2,\dotsc, -e_n)$, $0 \leq t \leq 2\pi$.

Since $\pmb{e} \ne \pmb{\overline{e}}$, the path $\delta(t)$ projects to a non-zero element $\alpha \in \pi_1 (V_{n,n-1})$,  which in turn is represented by the canonical Givens loop as in the above diagram $(G)$. The geometric Givens loop $\gamma^s(\theta)$ acts on the {\it decks} of the double covering. But $\gamma^s(\theta)$ acting on each $e_j$ gives just $e_j$, hence the action of the class $\alpha = \rho_n \left(\delta(t)\right) = \left[\gamma^s\right]$ on the ``sheets'' of $O(n)$ is trivial. This contradicts the observation that $\alpha \ne 0$. See \cite{Munkres}, Theorem 54.4, the "Fundamental Theorem of Lifting" for a covering mapping. \hfill $\blacksquare$

\bigskip
Note that in the above argument that $\gamma^s$ consists of two $n$-vectors $\cos 2\theta\,\, e_2 -\sin 2\theta\,\, e_3$ and $-\sin 2\theta\,\, e_2 + \cos 2\theta\,\, e_3$, taking $n=3$, which lie in the $y,z$ plane.

If the initial frame in ``sheet 1'' is $e_1, e_2, e_3$, then $e_1$ cannot move continuously to $-e_1$ under the action of $\gamma^s$. This reasoning applies for all $n\geq 3$, since in $O(n)$, the ``initial vectors'' $e_1, e_4, e_5,\dotsc, e_n$ {\it all remain invariant} under the ``deck transformation'' induced by $\gamma^s(\theta)$, see \cite{Sieradski}.

\head Davenport's Vector and a Perturbed Rotation\endhead

   In aeronautical guidance and related areas \cite{Shuster \& Oh}, the {\it Wahba problem} seeks to determine a rotation $A$ of $\RR^3$, that takes a collection of three vectors to itself, where the vectors are expressed via two distinct coordinate systems, say an observation frame $\{\text{Sun, Moon, Sirius}\}$ and a {\it natural} (reference) frame given by the body of the spacecraft. The given data consists of $N$ ``observation'' vectors $\{w_i\}$ and $N$ ``reference'' vectors $\{v_i\}$ so that ideally $Av_i = w_i$. But in practice the problem cannot be solved exactly since all vector estimates must be assumed to be corrupted by measurement error. A minimization formulation is often given as
   $$L(A) = \frac{1}{2} \sum_{i=1}^N \left|w_i - Av_i\right|^2\qquad,$$
   where we have suppressed optional ``measurement weights''.

   In earlier days, full eigenvector solutions were not considered feasible under high measurement rates or limited capacity of space-borne calculators. In the initial problem write-up \cite{Wahba}, Farrell and Stuelpnagel proposed a solution through the Polar decomposition as in Section 1. Specifically, the ``loss function'' $L(A)$ can be minimized through the {\it maximization} of
   $$F(A) = \text{trace} \left(A^T W V^T\right)\qquad,$$
   where $W$ means the rectangular matrix of ``observations'', and $V$ represents the rectangular matrix of ``body reference'' (both representing approximations of the same celestial directions).

   Considering now the square matrix $G = W V^T$ in Polar form as $G = UP$, we derive $F(A) = \text{trace} \left(A^T U P\right)$, or
   $$F(A) = \text{trace} \left(A^T U N^T D N\right) = \sum_{j=1}^3 d_j x_{jj}\qquad.$$
   Here we used the {\it real spectral} theorem to transform positive-definite $P$ (as arises generally) into $NPN^T = D$, where $N_{3 \times 3}$ is orthogonal and $D$ is real diagonal. Thus it is seen that maximization of $F(A)$ takes place when $X = I_{3 \times 3}$ in practice.

   Actually, for an optimal answer $A_0$, this $3 \times 3$ matrix should simply be $U$, the ``orthogonal part'' of $G_{3 \times3} = WV^T$. The close connection between ``polar form'' and $SVD$ also allows for a well known solution of Wahba's problem through singular-value decomposition of $WV^T$ as above, see \cite{Markley \& Montari}.

   Of course both $SVD$ algorithm technology and the {\it chips} that make it work have improved greatly since the 1960's, so rapid position updating for guidance is now routine.

   Instead of orthogonal $A$, one could deal with a unit quaternion vector $q$ such that $\Phi (q) = A$. Accordingly, P. Davenport, (see [Shuster], [Wahba]), showed that the loss function $L(A) = g(q)$ can be expressed as a quadratic form $g(q) = q^T Kq$, where $K$ is a $4 \times 4$ real symmetric matrix given in terms of $B = WV^T$ as above. This leads to the ``optimal'' quaternion $q_0$ (or $-q_0$) corresponding to the optimal Wahba $A_0$. In fact one may solve for
   $q_0$ directly, by finding the eigenvector for the largest eigenvalue of $K$. Iterations for solving
   $$K q = \lambda_{\max}\, q$$
   are natural in the context of the Spectral Theorem.

   We exhibit a specific case of the $K$-matrix that is useful toward our purposes.
I.Y. Bar-Itzhack shows that with {\it two} measurements, a given matrix which is {\it exactly orthogonal}
$$A = \left[\matrix
d_{11} & d_{12}& d_{13} \\
d_{21} & d_{22} & d_{23} \\
d_{31} & d_{32} & d_{33}
\endmatrix\right]$$
leads to the $4 \times 4$ matrix
$$ K_2 = \frac{1}{2}\left[\matrix
d_{11}-d_{22} & d_{21}+d_{12} & d_{31} & -d_{32} \\
d_{21}+d_{12} & d_{22}-d_{11} & d_{32} & d_{31} \\
d_{31} & d_{32} & -d_{11}-d_{22} & d_{12}-d_{21} \\
-d_{32} & d_{31} & d_{12}-d_{21} & d_{11}+d_{22}
\endmatrix\right]\qquad .$$
Note that $K_2$ is actually a function just of the {\it first two columns} of $A$. Hence $K_2$ is well-defined on our canonical ``path-component'' $V_{3,2}$, and so is the ``maximal'' eigen-vector of  $K_2$. What we have just given is a { \it more constructive} proof that $\Phi: \HH_1 \to V_{3,2}$ is a surjective mapping {\it without} using Brouwer's Invariance of Domain, see Appendix 0.

By contrast, in case an ``initial'' matrix $A$, or one which arises in an intermediate calculation, is {\it not} precisely orthogonal, one may ``correct'' the matrix by using quaternion methods. This approach emerged some years subsequent to the Davenport ``Lagrange multiplier'' or ``spectral'' solution of Wahba's problem. Major contributors were \cite{Klumpp}, \cite{Shepperd} and \cite{Landis Markley}. Amazingly, this concrete ``orthogonalization'' is given by a matrix that is {\it rational} in the original entries, with terms of degree at most {\it two} in both numerator and denominator.

Using Octave notation, starting with the symbolic matrix $S$ in Figure X we have
$$\left[\matrix
S11, S12, S13\\
S21, S22, S23\\
S31, S32, S33
\endmatrix
\right]
\tag{A}$$
and the generalized Landis matrix
$$\left[\matrix
S11 + S22 + S33 + 1,           S32 - S23,           S13 - S31,           S21 - S12\\
S32 - S23, S11 - S22 - S33 + 1,           S12 + S21,           S13 + S31\\
S13 - S31,           S12 + S21, S22 - S11 - S33 + 1,           S23 + S32\\
S21 - S12,           S13 + S31,           S23 + S32, S33 - S22 - S11 + 1
\endmatrix
\right]
\tag{B}$$
Working only with this first column called $ww1$, we may employ the ``conversion'' by $\Phi(q)$ giving matrix
$$D_w = \left[\matrix
\left( S11 + S22 + S33 +1 \right)^2 \\
& \ddots\endmatrix \right] \qquad. \tag{C}$$
$$D_w = {4 \times 4 \text{ matrix in Figure X\qquad.}}$$
To normalize this matrix, we must divide by $\ell_{11}^2 + \ell_{21}^2 + \ell_{31}^2 + \ell_{41}^2$, the norm=square of the column $ww1$. Furthermore, $ww1$ when {\it not} normalized, must be an eigenvector of the Bar-Itzhack matrix. The example $D_S$ [Itzhack Diary One] leads to $R_S$ the $K_2$ matrix.

We observe that when we apply ``landis'' to $D_S$, the first column is now called $LDq$. The adjusted quaternion $q' = [q_2, q_3, q_4, -q_1]$ is defined and ``$R_S$'' is applied so that $q'$ is seen to be an approximate eigen-vector of $R_S$, which equals the $4 \times 4$ Bar-Itzhack matrix of the rotations $D_S$, called $K_2$.

The symbolic verification that the first column $ww1$ of the Landis matrix for $S$ gives $\lambda \cdot ww1$, from $R_S \cdot ww1$, where $R_S$ is the Bar-Itzhack $K_2$ coming from symbolic $S$, is tendered to consideration of the Reader.

Now if a ``rotation matrix'' $D_S$ is precisely orthogonal, then the $4 \times 4$ Landis matrix construction leads to an {\it exact} quaternion direction $(q_1, q_2, q_3, q_4)^T$ which reverts to the same $3 \times 3$ matrix $D_S$ through the ``conversion'' $\Phi (\pmb{q})$. In fact this {\it exact case} gives, see \cite{Landis Markley},
$$\align
L(A) &= \text{Landis}(A)  =  \\
&\left[\matrix
A_{11} + A_{22}+A_{33}+1 & A_{32}-A_{23} & A_{13}-A_{31}  & A_{21}-A_{12} \\
A_{32}-A_{23} & A_{11}-A_{22}-A_{33}+1 & A_{12}+A_{21} & A_{13}+A_{31} \\
A_{13}-A_{31} & A_{12}+A_{21} & A_{22}-A_{11}-A_{33}+1 & A_{23}+A_{32}\\
A_{21}-A_{12} & A_{13}+A_{31} & A_{23}+A_{32} & A_{33}-A_{22}-A_{11}+1
\endmatrix\right] \\
& = 4 \cdot \left[\matrix
q_1^2 & q_1q_2 & q_1q_3 & q_1q_4 \\
q_2q_1 & q_2^2 & q_2q_3 & q_2q_4 \\
q_3q_1 & q_3q_2 & q_3^2 &q_3q_4 \\
q_4q_1 & q_4q_2 & q_4q_3 & q_4^2
\endmatrix\right] \qquad.\\
\endalign
$$
In case $A_{3 \times 3}$ is a {\it perturbation} of an orthogonal $\hat{A}$, we do not have a {\it unit} quaternion to make these formulas hold. We may use the Landis matrix $L(A)$, which possesses all four columns {\it parallel}, to create such a quaternion $q$. Let $v^{(1)}$ be the first column of $L(A)$ and set $q^{(1)} = \dfrac{v^{(1)}}{\|v^{(1)}\|}$, so we {\it divide} by a {\it square root}, which for aeronautics is an undesirable aspect of Klumpp's method, Shepperd's method, Markley's method and so forth. To obtain a good orthogonal approximation by $A$, we may avoid the need for any ``roots'' by noting that {\it every entry} of $\Phi(A)$ is {\it quadratic homogeneous}.

Thus we gain the rational expression $D_w = s \circ 3 (ww1)$ divided by $\Gamma = \text{nrmsq} (ww1)$ where we are treating the symbolic matrix $S$.

The assiduous Reader is invited to verify through symbolic algebra (polynomial division) that $D_{w}/\Gamma = S$, where $D_{w} = s o 3 (ww1)$, and $ww1$ is the first column of the $4 \times 4$ matrix landis$(S)$. That is, we may divide by $\Gamma$, the norm-square of the Landis column, to obtain a precise rotation $S$.
Here the ``assumptions'' inherent to ``symbolic orthogonal'' $S$ include $S^{(i)} \cdot S^{(j)} = \delta_{ij}$ for the columns $S^{(i)}$, $i =1,2,3$ and possibly similar formulas for the rows of $S$.

It is easier to work with the formula above for Land$(S)$ as $4 \cdot Q$, where $[Q]_{i,j} = q_i q_j$. We have all the monomial terms needed in the ``conversion'' $\Phi(q)$ to arrive back at the originating rotation, which should look exactly like
$$S = \left.\matrix
S_{11} & S_{12} & S_{13} \\
S_{21} & S_{22} & S_{23} \\
S_{31} & S_{32} & S_{33} &\quad .
\endmatrix\right.$$
This ``invertible'' linear tranformation, in terms of atoms $\{q_i q_j\}$, is somewhat unusual in that the Landis matrix $\text{Land}(S)_{4 \times 4}$ has column rank $=1$.

\newpage

{{\tt

Itzhack Diary One

\bigskip

Ds =

   -0.5450    0.7970    0.2600

    0.7330    0.6030   -0.3130

   -0.4070    0.0210   -0.9130

Rs =

   -0.5740    0.7650   -0.2035   -0.0105

    0.7650    0.5740    0.0105   -0.2035

   -0.2035    0.0105   -0.0290    0.0320

   -0.0105   -0.2035    0.0320    0.0290

>> LD=landis(Ds)

LD =

    0.1450    0.3340    0.6670   -0.0640

    0.3340    0.7650    1.5300   -0.1470

    0.6670    1.5300    3.0610   -0.2920

   -0.0640   -0.1470   -0.2920    0.0290

>> LD1=LD(1:4,1)

    0.1450

    0.3340

    0.6670

   -0.0640

>> LDq=LD1/sqrt(nrmsq(LD1))

LDq =

    0.1901

    0.4380

    0.8746

   -0.0839

>> LDr=LDq([2 3 4 1])

    0.4380
    0.8746
   -0.0839
    0.1901

>> LDr(4)=-LDr(4)

LDr =

    0.4380
    0.8746
   -0.0839
   -0.1901

>> Rs*LDr

ans =

    0.4368
    0.8749
   -0.0836
   -0.1908

>> LamLDr=ans./LDr

LamLDr =

    0.9973
    1.0003
    0.9961
    1.0034

}

\medskip
\centerline{Eigen-vector of the Bar-Itzhack $K\_2$ Operator}


\newpage

{{\tt

\bigskip

>> S=sym('S)

S =

[ S11, S12, S13]

[ S21, S22, S23]

[ S31, S32, S33]

>> w1=1+S(1,1)+S(2,2)+S(3,3)

w1 =

S11 + S22 + S33 + 1

>> assumeAlso(S,'real')

>> W

 W =

 [ S11 + S22 + S33 + 1,           S32 - S23,           S13 - S31,           S21 - S12]

[           S32 - S23, S11 - S22 - S33 + 1,           S12 + S21,           S13 + S31]

[           S13 - S31,           S12 + S21, S22 - S11 - S33 + 1,           S23 + S32]

[           S21 - S12,           S13 + S31,           S23 + S32, S33 - S22 - S11 + 1]

 >> ww1=W(1:4,1)

 ww1 =

  S11 + S22 + S33 + 1

           S32 - S23

           S13 - S31

           S21 - S12

 >> Dw=so3(ww1)

 Dw =

 [ (S11 + S22 + S33 + 1)\^{}2 - (S12 - S21)\^{}2 - (S13 - S31)\^{}2 + (S23 - S32)\^{}2,         2*(S12 - S21)*(S11 + S22 + S33 + 1) - 2*(S13 - S31)*(S23 - S32),         2*(S13 - S31)*(S11 + S22 + S33 + 1) + 2*(S12 - S21)*(S23 - S32)]

[       - 2*(S12 - S21)*(S11 + S22 + S33 + 1) - 2*(S13 - S31)*(S23 - S32), (S11 + S22 + S33 + 1)\^{}2 - (S12 - S21)\^{}2 + (S13 - S31)\^{}2 - (S23 - S32)\^{}2,         2*(S23 - S32)*(S11 + S22 + S33 + 1) - 2*(S12 - S21)*(S13 - S31)]

[         2*(S12 - S21)*(S23 - S32) - 2*(S13 - S31)*(S11 + S22 + S33 + 1),       - 2*(S23 - S32)*(S11 + S22 + S33 + 1) - 2*(S12 - S21)*(S13 - S31), (S11 + S22 + S33 + 1)\^{}2 + (S12 - S21)\^{}2 - (S13 - S31)\^{}2 - (S23 - S32)\^{}2]

}

\medskip
\centerline{Figure X}


\newpage

\head Conclusion\endhead

We examine the path connection of the compact Stiefel space of $n$-frames in\break $n$-space $\RR^n$ by means of well-known elementary facts from topology. The foundation for work in higher dimensions $n$, occurs in dimension 3, where the {\it group structure} upon $\HH_1 \simeq \SS^3$, the group of unit quaternions, comes into play.

One observation is that a {\it perturbed} rotation can be straightened out to a nearby rotation {\it minimal} in the sense of the angular adjustment required. Such formulas are indispensable in the practice of aeronautical guidance. The related ``Procrustes'' problem seeks such a closest rotation in the Frobenius norm.
In 3-space the Procrustes question may be solved through an eigen-vector formulation proposed by \mbox{P\.~Davenport}.  Concrete expressions due to I.Y. Bar-Itzhack have the theoretical impact that $\Phi : \SS^3 \to SO (3)$ is a continuous surjection. The fact that the Davenport eigen-vector can always be found, bears on our original ``path-component'' problem for frames in 3-space: the fact that the frames are classifiable into left-hand and right-hand path-components.

\newpage



\newpage

\head
Appendix 0
\endhead

\subhead
Domain Invariance on Manifolds
\endsubhead

The following expresses roughly that within the category of m-manifolds, a locally one-to-one mapping is also a covering map: it is locally a homeomorphism.

\bigskip
\noindent
{\bf Proposition 0-1}\quad
Let $M, N$ be compact (Hausdorff) m-manifolds with $\tau : M \to N$ locally an injection.
If $N$ is path-connected, $\tau$ must be (``globally'') surjective, that is, $\tau (M) = N$.

\bigskip
\noindent
{\it Proof.} \quad
Under the compactness hypothesis, $\tau (M)$ must be closed in $N$, \cite{Dugundji}.

Let $U \subset M$ be an open subset homeomorphic to $\RR^m$, so that $\tau|_U$ induces a {\it one-to-one} mapping onto $\tau (U) = W \subset N$\,. As in \cite{Deimling} p. 23, this restricted mapping must be {\it open}, hence $\tau (M)$ is open {\it and} closed in $N$ so must cover surjectively at least one of the path-components of $N$. This proof uses the Brouwer degree of an {\it ``odd''} (antipode-preserving) map. Alternatively, any {\it odd} self-mapping of a sphere is {\it essential}. This goes back to versions of B$\acute{e}$zout's theorem on real varieties inspired by [Behrend], \cite{Kapferer} and \cite{Pfister}. The connection between B$\acute{e}$zout's theorem and Invariance of Domain is outlined in \cite{Sjogren, Homogeneous}.
\hfill $\blacksquare$

\newpage

\head
Appendix 1
\endhead

\subhead
The Symbolic Landis Orthogonalization
\endsubhead

If an abstract ``rotation'' be given as

$$S = \left[\matrix
S11, S12, S13 \\
S21, S22, S23 \\
S31, S32, S33 \\
\endmatrix\right]\qquad,$$
then a $4\times 4$ matrix Landis(S), written

$$Q_c=\text{landis}(S)$$

$$Q_c  =
\left[\matrix
S11 + S22 + S33 + 1, S32 - S23, S13 - S31, S21 - S12 \\
S32 - S23, S11 - S22 - S33 + 1, S12 + S21, S13 + S31 \\
S13 - S31, S12 + S21, S22 - S11 - S33 + 1, S23 + S32 \\
S21 - S12, S13 + S31, S23 + S32, S33 - S22 - S11 + 1
\endmatrix\right]$$
can in turn be ``converted'' according to $\Phi(w1)$ where $w1$ is the first {\it column} of $Q_c$ (or use another  column).

The resulting symbolic matrix, algebraically expanded, is given as $\mbox{D}_{\mbox{S}}$, before we divide by the scalar polynomial homogeneous of degree two, namely
$$\parallel w1 \parallel^2 {:} {\,=\,} \text{denom}\qquad .$$

$\mbox{D}_{\mbox{S}}$ =

[ 2*S11 + 2*S22 + 2*S33 + 2*S11*S22 + 2*S12*S21 + 2*S11*S33 + 2*S13*S31 + 2*S22*S33 - 2*S23*S32 + 4*S22\^{}2 + 4*S23\^{}2 + 4*S32\^{}2 + 4*S33\^{}2 - 4,       2*S12 - 2*S21 + 4*S12*S22 - 4*S21*S22 + 2*S12*S33 + 2*S13*S32 - 2*S21*S33 + 2*S23*S31 - 4*S31*S32,       2*S13 - 2*S31 + 2*S12*S23 + 2*S13*S22 - 4*S21*S23 + 4*S13*S33 + 2*S21*S32 - 2*S22*S31 - 4*S31*S33]

[                                         2*S21 - 2*S12 - 4*S12*S22 - 4*S13*S23 + 4*S21*S22 - 2*S12*S33 + 2*S13*S32 + 2*S21*S33 + 2*S23*S31, 2*S11 + 2*S22 + 2*S33 + 2*S11*S22 + 2*S12*S21 + 2*S11*S33 - 2*S13*S31 + 2*S22*S33 + 2*S23*S32 + 4*S22\^{}2,                   2*S23 - 2*S32 + 2*S11*S23 + 2*S13*S21 - 2*S11*S32 + 2*S12*S31 + 4*S22*S23 + 4*S23*S33]

[                                         2*S31 - 2*S13 + 2*S12*S23 - 2*S13*S22 - 4*S12*S32 - 4*S13*S33 + 2*S21*S32 + 2*S22*S31 + 4*S31*S33,                   2*S32 - 2*S23 - 2*S11*S23 + 2*S13*S21 + 2*S11*S32 + 2*S12*S31 + 4*S22*S32 + 4*S32*S33, 2*S11 + 2*S22 + 2*S33 + 2*S11*S22 - 2*S12*S21 + 2*S11*S33 + 2*S13*S31 + 2*S22*S33 + 2*S23*S32 + 4*S33\^{}2]\qquad

\bigskip
But if we write $w1 = 4 q_1 [q_1, q_2, q_3, q_4]$
we obtain re-Markley

$$\text{denom} = 16 q_1^2$$
and also $4 q_1^2 = 1 + S_{11} + S_{22} + S_{33}.$
Thus $\parallel w1 \parallel^2 = \text{nrmsq} (w1) = 4\,(1 + S_{11} + S_{22} + S_{33}).$

These identities relate to the fact that for the secular polynomial is one variable $X$,
$$\text{secul} (S) = X^3 - \text{trace} (S) \cdot X^2 + \text{principal minors of }(S) \cdot X - \det (S)\qquad.$$

In this Appendix we are not so pedantic in handling all $3 \times 3$ orthogonal matrices at once; the formulas in the case of a reflection-rotation differ from those of a rotation, by a sign change.

As an example starting from the rotation

$$t =
\left.\matrix
-0.6651  &  0.7463  &  -0.0256\\
-0.7395  & -0.6631  & -0.1162\\
-0.1037  & -0.0583  &  0.9929
\endmatrix\right.$$
we compute Landis(t) as a symmetric matrix of {\it rank one}.

$$wt=\text{landis}(t)$$

$$wt =
\left.\matrix
0.6647  &  0.0578  &  0.0781  & -1.4858\\
0.0578  &  0.0050  &  0.0068  & -0.1293\\
0.0781  &  0.0068  &  0.0092  & -0.1745\\
-1.4858 & -0.1293  & -0.1745  &  3.3211
\endmatrix\right.\qquad .$$

From the first column of Landis(t), there are two ways to compute $\parallel w1 \parallel^2$

$$w1=wt(1:4,1)$$

$$w1 =
\left.\matrix
& 0.6647\\
& 0.0578\\
& 0.0781\\
& -1.4858
\endmatrix\right.\qquad,$$
giving $\text{nrt} = \text{nrmsq(w1)} = 2.6588$.

The second one is

$$4 \cdot (1 + t_{11} + t_{22} + t_{33})\qquad,$$
which yields 2.6588 as well.

Either way, the original rotation $t$ is recovered by the conversion formula. Given an {\it approximate} rotation $\tau$, one may experiment with trying the apparently {\it linear} denominator $4 \cdot (1 + \tau_{11} + \tau_{22} + \tau_{33})$ to obtain a reasonable approximation to the ``closest orthogonal version of $\tau$'', up to machine precision.  One may re-scale the matrix rows toward unit norm.


\newpage

\head
Appendix 2
\endhead

\subhead
Correspondence with John Stallings (1996)
\endsubhead

(to J.R. Stallings)     The bulk of this note will involve a slight application of topology to a linear algebra question, as if it applies. It was originally typed as an email, but since then some of the technical symbolism got typeset.

\bigskip
{\bf Remark.} These notes indicate the initial idea of establishing the parity of rotation-frames without needing the ``determinant'', that should be defined and continuous on an open set of matrices. Even a complete text such as Greub, "Linear Algebra" does not give details of the classical proof. It is hoped that this Appendix could provide motivation for the body of this article, which tries to address these intuitions with greater rigor.

\bigskip

OK let us look at my math problem.  All I want at this point is your sense that I am arguing plausibly here.

You are the right person too since you know about fiber spaces, plus you are an authority on exterior product constructions and the like in rings and modules.  But fear not, all we deal with is vectors in Euclidean space.  Probably all the topology needed is found in Massey's junior/senior textbook.

I am explaining “determinant” of a real matrix by the effect of the matrix on shapes and volumes in $R^N$.  At some point I decide that all volumes are positive, so then we introduce the concept of orientation of “frame” (say defining vectors of an $N$-polyhedron).

We have in $R^3$ that ($e1,e2,e3$) is rotatable into ($e3,e1,e2$) but not into ($e2,e1,e3$).  Rotations can be composed of planar rotations through an angle (Givens rotations).

By explicit construction in $R^N$ one can get than any orthogonal set ($f1, \cdots fN$) of vectors is rotatable either to

$E= (e1,e2, …. eN)$ or to  $EE=(e2,e1, e3, \cdots eN)$.  If fi is not a unit vector you can make it so simultaneously with the other through a family of scalings (homotopy to shrink or extend).

So under (rigid) rotation, the orthogonal “frames” consist of at most two classes.  We want to show that there are at least two classes.

Let a “frame” be a set of N N-vectors in $R^N$ .  We mostly are concern with “good” frames, that is, linearly independent frames.
Without any loss of generality, all vectors can be unit vectors as null vectors in a frame never enter the discussion (they might enter as linear combinations of frame vectors).

So we change the rotatability problem to a homotopy problem. If E were rotatable to EE, then there would be a homotopy leading from the good frame E to the good frame EE, taking values only in *good frames.  So if we prove that there is no such homotopy of E to EE within the space of good frames, they are not in the same rotation class either.  To go from a “good” frame to a an orthonormal frames, apply the QR algorithm.

You would think that dealing with rotations of orthogonal frames would be easier than homotopies of good frames, but in considering orthogonal frames, you are inevitably led back to a matrix formulation.  The temptation to prove this using continuity properties of the determinant function are irresistible, but this is not allowed in my context, since we are developing a foundation for the meaning of determinant.  Besides, this may show a (rare) application of elementary topology.  The standard way to show that two bases are homotopic if and only if their transformation matrix has positive determinant can be found on page 76 of “Linear Algebra” by W.H. Greub. Note however that Greub applies the intermediate value property to the determinant (as a function of quantity N vectors each in $R^N$), and invokes its continuity.  But he does not prove that continuity, at least in this book.

All right, now restricting to all unit vectors, we can categorize the space of all frames in $R^N$.  Each vector lies on a copy of $S^n$, so $BF = X S^n$, where $n=N-1$    ($N \times$, a Cartesian sum).  On the other hand, we may call the space of good frames,  $BG < BF$.  Here we let $N>1$ and handle trivial cases separately.

Let us build up the good frames as the total space of a fiber space.  Assume for the moment that we have done this, and have the topological space $F_n$ consisting of quantity $n=N-1$ linearly independent unit vectors {$v_i$}.  To get a “good frame” you need one more vector on $S^n$, subject to one constraint.  The given partial frame, or n-frame, in $F_n$, consists of quantity n unit vectors, and these generate a Great $S^{n-1}$ on the unit sphere $S^n$.  Pick u in $S^n$ which does *not lie on this Great Sphere.  Then $(v_1, \cdots v_n, u)$ constitutes a good frame.  Note that we have now constructed BG as a fibration with base
$F_n$ and fiber $W = S^n  \setminus S^{n-1}$, the n-sphere with its equator removed.  Clearly W has two topological components.  I should have called the $F_i$,  $G_n$ and $G_i$ for “good” but I don't want to change notation now.  $F_n$ is the base of the fibration.

The question is, does BG have more than one topological component, or is it connected?  If it is not connected, we saw by reducing the frames as above (doing some work probably with Gram-Schmidt orthogonalization), that there are exactly two components.  We observe that W does have two components, an upper hemisphere H+ and lower one H- .  The question is whether there might be a path from some point on H+ to a point on H-, fixing the base point c in $F_n$ but letting the path wander through $F_n$.  In that case any two frames are path equivalent and we can’t prove they are not rotation equivalent; in fact by Gram-Schmidt (the homotopy version using Givens rotations) you could also prove that any two frames are rotation equivalent.

We could have defined $F_i$ as space of quantity “i” independent unit vectors.  To get $F_{i+1}$ we started with $v_1 \cdots v_i$ and look for a unit vector not in their span, in $S^n$.  In other words the quantity “i” vectors generate a “principal” $S^{i-1}$ . [By “principal sphere” you would mean, “the equator of an equator of an equator” etc]    Choose any $v_{i+1}$ which is not in this principal sphere. Then $S^{i-1}$ is generated as the intersection of the plane $L = \text{span}\{v_1, \cdots, v_n\}$ and $S^n$.

So we build up $F_j$ as a sequence of fibrations whose fibers are $S^n$ with some lower-dimensional equator removed.  Remember that $n>=1$.  I suppose that this construction is very standard, was done long ago in some thesis, and the results are way beyond this trivial write-up.  But the application to simple linear algebra may be less well-known.

Now $F_1$ is just $S^n$ minus the ${-1}$ sphere, which is the empty set.

We claim that for $j<n$, all $F_j$ are connected and simply connected, and that $F_n$ is connected.

In fact we note that the fiber used in forming $F_{i+1}$ is $S^n  \setminus S^{i-1}$ .  In the case of forming $F_n$, the fiber was $S^{n} \setminus S^{n-2} = Y$ .  We know that $Y$ is connected but $\pi_1(Y) = Z$  generated by a circle $D$ constructed as follows.

Given $R= S^{n-2}$, it is the equator of some equator $Q=S^{n-1}$.  Thus choose the poles of R on Q, a+ and a$-$ (antipodal on Q and $S^n$)  On the other hand Q has unique poles b+ and b$-$ on $S^n$.  the two antipodal pairs are distinct and generate a principal circle $D = sp(a+, b+)$. The circle D itself is independent of the choices of Q etc., and depends only on R.

Now we look at page 377 of “Spanier” and note the long homotopy sequence of a (weak) fibration.

The important part is\vskip-18pt
$$\pi_1(fiber) -> \pi_1 (total) -> \pi_1(base) ->\pi_0(fiber) ->\pi_0(total)->\pi_0(base)\qquad .$$

\vskip-6pt

It is not obvious how those $\pi_0$ maps can be regarded as homomorphisms, but we won’t worry about that (explained on pg 371).  This mathematics is not too involved, as we notice that $\pi_2$ terms don’t ever really enter into the deliberations. The other fibration approach is based on $\pi_2$ of the spheres.

At each successive stage $F_i$ to $F_{i+1}$ the old total space becomes the new base space, so the total space has $\pi_1$ and $\pi_0$ equal to 0 at least up until $\pi_1$(fiber) not eq 0 which occurs when base is $F_{n-1}$ and fiber is the space Y, the N sphere minus an equator of an equator.

In that case we get that $F_n$ is still connected, but its $\pi_1$ is a quotient of $\pi_1$ of the fiber Y, which we saw was the additive integers Z, generated by that circle D (regarded as a path or loop with suitable base point).  If we want to be specific, call that fundamental quotient group, Zg.

Thus the whole fundamental group of $F_n$ arises from that fiber Y, and is generated by D, so if $\phi:F_n \to  F_{n-1}$ is the fibration map, we can take $\phi^{-1} (x_0) =Y$, sitting in $F_n$.

Summarizing we got $W < F_N \to F_n$, with fibration map $\psi$ where $W = \psi^{-1}(t_0)$ has two contractible components.  $F_n$ is connected and has $\pi_1$ equal to a quotient of the group of integers Z, carried by a loop in $Y < S^n \setminus WS^{n-2}$.

Start with $t_0 ={w_1, \cdots, w_{n-1}, u_0}$, then W is our fiber.  We keep $w_1 \cdots w_{n-1}$ fixed, these are quite arbitrary and can be chosen as an orthonormal $n-1$ set of unit vectors if desired.  The last vector u will vary. We have D as a principal circle in $Y = S^n \setminus WS^{n-2}$ where the sphere we removed was generated by the ${w_j}$.  Choose $b_0$ at an angle $\pi/2$ radians along that circle. Now define a smooth path gamma(s) of unit speed along D, with initial value $\gamma(0)=u_0$ and initial direction, in the direction of $b_0$.  The path gamma will terminate at $\gamma(1)= u_0$ again, and is thus a loop representing a generator of $\pi_1$ of Y and thus $\pi_1$ of $F_n$ (the l.i. sets of quantity n unit vectors).  For all we know. $\pi_1(F_n)$ might be the trivial group.

Now consider that all spheres $S^n$ have the same coordinatization.  For each speed parameter s, we have $\gamma(s) = u$, which defines a fiber $W_u$  (since we took all ${w_j}$ as fixed).  Furthermore a vector b is associated to $u$ by $u \to b$ by the following:

1.  b is on the principal circle called D

2.  the polar distance from u to b stays the same (90 degrees)

Note that while $u_0$ does not live in the fiber $W_{u_0} ~= S^n \setminus WS^{n-1}$, $b_0$ does so, and is in one of the connected components of $W$, say $H+$.  Similarly b is well-defined in $W_u$ corresponding to the speed parameter s.

What is your opinion about this loop construction?  I was a little worried that argument that the component of b remains fixed under the loop action, depends in some way on the orientation of the coordinate frame of $S^n$, interfering with the actual result we are trying to prove, and thus that the reasoning we are using is somehow “circular”.  Ha !  Let me know if this seems all right to you..

But that is what I am arguing.  Think about 2 vectors in the plane, leading to the 1-sphere in this discussion …

At $\gamma(1)= u_0$ again, the corresponding b must be $b_0$ again, since its angular distance to u was always $\pi/2$.  Thus we have a path lambda: $[0,1] \to BG$, such that $\lambda(0)=\lambda(1)= b_0$.  Thus $\lambda(1)$ is in the same connected component as $\lambda(0)$ so the induced deck transformation is “trivial” (the identity).  The path lambda projects to $F_n$ and in fact to the fiber Y, giving the same fundamental class as gamma, shifting the base point by $\pi/2$.  By the Fundamental Theorem of Covering Spaces, now any generator of $\pi_1 (F_n)$ induces a trivial deck transformation on the covering space whose fiber is the pair of components ${H+, H-}$ over every point of $F_n$.

Therefore there can never be a path in the total space BG of that fibration, leading from a point in the hemisphere H+ to the hemisphere H-.  Thus the fact that the fiber Y is not connected implies that BG also has more than one component, which we were trying to prove.  In fact BG always has two connected components represented by some “even” frame such as $e1,e2,e3 \cdots$ and an “odd” frame such as $e2, e1, e3$ respectively.

\newpage

\head
Appendix 3
\endhead

\subhead
Octave Matrix and Rotation Routines
\endsubhead

{{\tt

\bigskip

\bigskip

function W=landis(S)

W=[];

w=[];

if size(S)==[3 3]

    w1=1+S(1,1)+S(2,2)+S(3,3);

    w2=S(3,2)-S(2,3);

    w3=S(1,3)-S(3,1);

    w4=S(2,1)-S(1,2);

    w=[w1 w2 w3 w4];

    W=[W w'];

    w2=1+S(1,1)-S(2,2)-S(3,3);

    w1=S(3,2)-S(2,3);

    w4=S(1,3)+S(3,1);

    w3=S(2,1)+S(1,2);

    w=[w1 w2 w3 w4];

    W=[W w'];

    w3=1-S(1,1)+S(2,2)-S(3,3);

    w4=S(3,2)+S(2,3);

    w1=S(1,3)-S(3,1);

    w2=S(2,1)+S(1,2);

    w=[w1 w2 w3 w4];

    W=[W w'];

    w4=1-S(1,1)-S(2,2)+S(3,3);

    w3=S(3,2)+S(2,3);

    w2=S(1,3)+S(3,1);

    w1=S(2,1)-S(1,2);

    w=[w1 w2 w3 w4];

    W=[W w'];

else W=0

end

\bigskip

\bigskip

function M=itzhak(S)

M=[];

m=[];

if size(S)==[3 3]

    m1=S(1,1)-S(2,2);

    m2=S(2,1)+S(1,2);

    m3=S(3,1);

    m4=-S(3,2);

    m=[m1 m2 m3 m4];

    M=[M m'];

    m1=S(2,1)+S(1,2);

    m2= -S(1,1)+S(2,2);

    m3=S(3,2);

    m4=S(3,1);

    m=[m1 m2 m3 m4];

    M=[M m'];

    m1=S(3,1);

    m2=S(3,2);

    m3=-S(1,1)-S(2,2);

    m4= -S(2,1)+S(1,2);

    m=[m1 m2 m3 m4];

    M=[M m'];

    m1=-S(3,2);

    m2=S(3,1);

    m3=S(1,2)-S(2,1);

    m4=S(1,1)+S(2,2);

    m=[m1 m2 m3 m4];

    M=[M m'];

    M=M/2;

else M=0;

end

\bigskip

\bigskip

function spin=so2(v)

spin=sym(eye(2));

spin(1,1)=v(1);

spin(1,2)=-v(2);

spin(2,1)=v(2);

spin(2,2)=v(1);

\bigskip

\bigskip

function spun=so3(v)

spun=sym(eye(3));

x=v(1);

y=v(2);

z=v(3);

w=v(4);

spun(1,1)=x*x+y*y-z*z-w*w;

spun(1,2)=2*(y*z-x*w);

spun(1,3)=2*(y*w+x*z);

spun(2,1)=2*(y*z+x*w);

spun(2,2)=x*x-y*y+z*z-w*w;

spun(2,3)=2*(z*w-x*y);

spun(3,1)=2*(y*w-x*z);

spun(3,2)=2*(z*w+x*y);

spun(3,3)=x*x-y*y-z*z+w*w;

\bigskip

\bigskip

function ns = nrmsq(v)

v1=v(1);

v2=v(2);

v3=v(3);

v4=v(4);

ns= v1*v1+v2*v2+v3*v3+v4*v4;

\bigskip

\bigskip

function [f g h k] = qprod(v,w)

v1=v(1);

v2=v(2);

v3=v(3);

v4=v(4);

w1=w(1);

w2=w(2);

w3=w(3);

w4=w(4);

f= v1*w1-v2*w2-v3*w3-v4*w4;

g=v1*w2+v2*w1+v3*w4-v4*w3;

h=v1*w3-v2*w4+v3*w1+v4*w2;

k=v1*w4+v2*w3-v3*w2+v4*w1;

\bigskip

\bigskip

function cv = conjq(v)

cv=sym([1 1 1 1]);

v1=v(1);

v2=v(2);

v3=v(3);

v4=v(4);

cv(1)=v1;

cv(2)=-v2;

cv(3)=-v3;

cv(4)=-v4;

\bigskip

\bigskip

function [f g] = cprod(v,w)

a=v(1);

b=v(2);

c=w(1);

d=w(2);

f= a*c-b*d;

g=a*d +b*c;

}

\end{document}

\end